\documentclass[12pt]{article}
\usepackage{amssymb}
\makeatletter

\@addtoreset{equation}{section}

\makeatother

\newtheorem{thm}[equation]{Theorem}
\newtheorem{lemma}[equation]{Lemma}

\renewcommand\thefigure{\thesection.\@arabic\c@figure}

\renewcommand\thetable{\thesection.\@arabic\c@table}

\def\reff#1{(\ref{#1})}
\def\sobre#1#2{\lower 1ex \hbox{ $#1 \atop #2 $ } }
\def\supp{{\rm Supp}\,}
\def\basis{{\rm Basis}\,}
\def\life{{\rm Life}\,}
\def\birth{{\rm Birth}\,}
\def\death{{\rm Death}\,}
\def\tl{{\rm TL}\,}
\def\ti{{\rm TI}\,}
\def\sw{{\rm SW}\,}

\begin{document}

\def\E{{\mathbb E}}
\def\P{{\mathbb P}}
\def\R{{\mathbb R}}
\def\Z{{\mathbb Z}}
\def\V{{\mathbb V}}
\def\N{{\mathbb N}}
\def\NN{{\bf N}}
\def\X{{\cal X}}
\def\Y{{\bf Y}}
\def\G{{\cal G}}
\def\bS{{\bf S}}
\def\T{{\cal T}}
\def\C{{\bf C}}
\def\D{{\bf D}}
\def\G{{\bf G}}
\def\U{{\bf U}}
\def\Ke{{\bf K}_{\zeta}}
\def\K{{\bf K}}
\def\H{{\bf H}}
\def\n{{\bf n}}
\def\b{{\bf b}}
\def\g{{\bf g}}
\def\mm{{m}}
\def\sqr{\vcenter{
         \hrule height.1mm
         \hbox{\vrule width.1mm height2.2mm\kern2.18mm\vrule width.1mm}
         \hrule height.1mm}}                  
\def\square{\ifmmode\sqr\else{$\sqr$}\fi}
\def\one{{\bf 1}\hskip-.5mm}
\def\liml{\lim_{L\to\infty}}
\def\given{\ \vert \ }
\def\Given{\ \Big\vert \ }
\def\be{{\beta}}
\def\de{{\delta}}
\def\la{{\lambda}}
\def\ga{{\gamma}}
\def\th{{\theta}}
\def\proof{\noindent{\bf Proof. }}
\def\rate{{e^{- \beta|\ga|}}}
\def\A{{\bf A}}
\def\B{{\bf B}}
\def\C{{\bf C}}
\def\D{{\bf D}}

\title{Loss network representation of Peierls contours}
\date{}
     
\author{Roberto Fern\'andez \\[-3mm] {\it \normalsize Consejo Nacional de
    Investigaciones Cient\'{\i}ficas y T\'ecnicas, Argentina}\\[-3mm] 
{\it \normalsize and Universidade de S\~{a}o Paulo}
\\ 
Pablo A. Ferrari\\[-3mm] {\it \normalsize Universidade de S\~{a}o Paulo} 
\\ 
Nancy L. Garcia \\[-3mm] {\it \normalsize Universidade Estadual de Campinas} 
}

\maketitle
\abstract{}\endabstract
\vskip-10mm
\baselineskip 17pt 
We present a probabilistic approach for the study of systems with exclusions,
in the regime traditionally studied via cluster-expansion methods.  In this
paper we focus on its application for the gases of Peierls contours found in
the study of the Ising model at low temperatures, but most of the results are
general.  We realize the equilibrium measure as the invariant measure of a
loss-network process whose existence is ensured by a subcriticality condition
of a dominant branching process.  In this regime, the approach yields, besides
existence and uniqueness of the measure, properties such as exponential space
convergence and mixing, and a central limit theorem.  The loss network
converges exponentially fast to the equilibrium measure, without metastable
traps.  This convergence is faster at low temperatures, where it leads to the
proof of an asymptotic Poisson distribution of contours. Our results on the
mixing properties of the measure are comparable to those obtained with
``duplicated-variables expansion'', used to treat systems with disorder and
coupled map lattices. It works in a larger region of validity than usual
cluster-expansion formalisms, and it is not tied to the analyticity of the
pressure.  In fact, it does not lead to any kind of expansion for the latter,
and the properties of the equilibrium measure are obtained without resorting
to combinatorial or complex analysis techniques.

\bigskip




\baselineskip 22pt
\section{Introduction}
In this paper we develop a probabilistic approach to the study of
the equilibrium measure of systems with exclusions ---like hard-core
gases, contours, polymers or animals--- in the low-density or
extreme-temperature regime.  This regime has traditionally been studied
via cluster-expansion methods, which relied either on 
sophisticated combinatorial estimations [Malyshev (1980), Seiler (1982),
Brydges (1984)] or on astute inductive hypotheses plus complex
analysis [Koteck\'y and Preiss (1986), Dobrushin (1996, 1996a)].  

In contrast, we realize the equilibrium measure as the invariant
measure of a loss network process that can be studied using standards
tools and notions from probabilistic models and processes.  Loss
networks, first introduced by Erlang in 1917, encompass a rather
general family of processes as discussed in Kelly (1991) and
references therein.  Technically we work with the so-called
\emph{fixed-routing} loss networks.  We build on ideas of
Ferrari and Garcia (1998) to show (Section \ref{S34}) that the
existence of the loss network can be related to the absence of
percolation in an oriented percolation process.  This condition also
yields other properties of the process and its invariant measure, like
uniqueness, convergence through sequences of finite volumes and mixing
properties (Theorem \ref{225}).  More precise results can be obtained
by resorting to a dominant multitype branching process (Section
\ref{SSS}).  Roughly speaking, the mean number of branches of this
process becomes the driving parameter: Subcriticality is a sufficient
condition for the construction to work.  Time and space rates of
convergence and mixing rates are explicitly obtained in terms of this
parameter (Theorem \ref{18a}).

The approach of this paper was already exploited in Fern\'andez,
Ferrari and Garcia (1998). Here we refine and complete the theory
presented there, extending their region of validity and including
proofs of exponential mixing and convergence to a Poisson process.

For concreteness, we analyze in this paper the gases of Peierls
contours used, for instance, for low-temperature studies of the Ising
model [Peierls (1936), Dobrushin (1965), Griffiths (1964)].  The
subcriticality condition of the corresponding branching process is
\begin{equation}
\sup_\gamma {1\over |\gamma|} \sum_{\theta : \theta\not\sim\gamma} |\theta|
\,w(\theta)\;<\;1\;,
\label{fal.1}
\end{equation}
where $|\gamma|$ indicates the length (perimeter, surface area) of a
contour $\gamma$, $w(\ga)$ is its weight  and ``$\sim$''
stands for the volume-exclusion (=non-intersection) condition.  

Condition \reff{fal.1} is considerably weaker than those obtained by
the most developed cluster-expansion approaches [Koteck\'y and Preiss
(1986), Dobrushin (1996, 1996a)] in which the factor $|\theta|$ in the
right-hand side is replaced by a function that grows exponentially
with $|\theta|$.  The weakening of the condition has a price: unlike
previous approaches, ours does not yield analyticity properties of the
expectations.  Condition \reff{fal.1} is similar to conditions
obtained in the study of systems for which nonanalyticity is known
[von Dreifus, Klein and Perez (1995)] or suspected [Bricmont and
Kupiainen (1996, 1997)].  Our results on space convergence and mixing
rates can basically be obtained with the ``duplicate-system''
expansions of Bricmont and Kupiainen (1996, 1997).

The novelty of our approach lies in the following features: First, it
offers a completely different framework for the study of hard-core
measures, based on well known stochastic processes.  This can
conceivably lead to new insights and stronger results. In particular
the condition \reff{fal.1} can potentially be weakened via
sub-criticality estimates obtained directly for the associated
oriented percolation process, without resorting to a dominating
branching process.  Similar improvements have been done for
the contact process and oriented percolation, see Liggett (1995), for
instance. Second, our construction involves a stochastic process (the
loss network) that converges exponentially fast to the sought-after
measure. This process is, in principle, easy to simulate and its
potential as a computational tool deserves to be explored
(Fern\'andez, Ferrari and Garcia, 1999).  Its rate of convergence to
the equilibrium measure increases as the temperature decreases and,
unlike spin-flip dynamics, it does not present meta-stable traps (any
contour lives an exponential time of mean one).  Finally, the
construction permits a rather straightforward proof of the asymptotic
Poisson distribution of contours at low temperature (Theorem
\ref{191}). It is not obvious to us how such a result can be obtained
through the standard statistical mechanical expansions.

The present approach does not lead to a series expansion for the
pressure (or free-energy density).  In particular, it does not yield
the ``surface-tension bounds'' that play such a crucial role in some
applications of cluster expansions [see eg.\ Zahradn\'{\i}ck (1984),
Borgs and Imbrie (1989)].  In fact, the approach is designed so to
\emph{bypass} expansions of this type.  It is a probabilistic approach
designed to answer probabilistic questions ---existence of
expectations, properties of correlation functions--- in a direct way,
without combinatorial or complex-analysis techniques.  It is not an
alternative to cluster expansions: it has a different regime of
validity and different aims.

While the ideas behind our results are natural and simple, their
formalization requires many intermediate technical results that may
obscure the development of our theory.  Let us, therefore, present a
sort of ``road map'' of the paper to guide the reader.  The main
results are presented, in a self-contained manner, in Section
\ref{SS2}.  The actual construction of the loss network is the subject
of Section \ref{S3}. We start with a reference ``free'' process
of Poissonian births and exponentially distributed deaths, with
respect to which the loss network is absolutely continuous. The
novelty of our approach resides in the fact that, rather than
independently generating birth and death-times, the lifetimes are
associated to each birth-time as a mark.  Hence, unlike comparable
constructions, each death-time has an associated birth-time and defines
an space-time ``cylinder'' representing the presence of a loss event
(contour).  This permits to comb the process either backwards or
forward in time with equal ease.  The ``time backwards'' point of view
leads to the notion of ``backwards oriented percolation'' which will
be our main conceptual and practical tool.  The idea is to construct
the loss network by erasing from the free process
those cylinders that conflict with pre-existing ones.  This can only
be done if the set of preexisting cylinders (the ``clan of
ancestors'') is finite with probability one.  This is precisely the
condition of absence of backwards oriented percolation.  The main
technical result of this Section \ref{S3} is the proof that this
construction yields precisely the loss network.  [This is contained in
the proof of Theorem \ref{20}]. It is apparent that the same
construction works for any other process absolutely continuous with
respect to a Poisson birth and death process, for instance point or
Boolean processes (Baddeley and van Lieshout, 1995; Kendall, 1997,
1998).

Once the process is so constructed, the course of action is clear.
First, we relate time and space mixing properties with the time and
space size of the percolation clan.  This is done in Section
\reff{SS4} [Theorem \ref{225}].  The only slightly involved part of
this section is related to the proof of the mixing properties.
Indeed, we resort to a ---standard but unusual--- coupling between
clans, together with a continuous-time construction (Section
\ref{ssfc}), to improve a little over results proved previously by the
method of ``duplicated variables'' (von Dreifus, Klein and Perez,
1995; Bricmont and Kupiainen, 1996).

In order to present quantitative estimates in terms of the parameters
of the problem, we follow the known technique (used, for instance by
Hall, 1985) of bounding percolation probabilities via a branching
process.  In particular, the size of the percolation clan is 
bounded by number of branches. This is done in Section \ref{SSS},
where the only (very simple!)  algebraic calculations of the paper are
presented [displays \reff{a16}, \reff{702} and \reff{802}].  As
expected, subcriticality of the branching process implies lack of
percolation and exponentially damped sizes of the percolation clan.
The main estimations are contained in Theorem \ref{thm:hall}.  Once
again, we present a slightly unusual continuous-time construction
(Section \ref{ctbp}) to improve one of the estimations, namely the
time-length of a clan [part (ii) of Theorem \ref{thm:hall}].  Readers
can opt instead for the more direct, but slightly weaker, estimate
presented in the first remark following the theorem.

The results of Theorem \ref{18a} are a direct consequence of the
estimates of Section \ref{SSS} applied to the percolation expressions
of Section \ref{SS4}, as explained in Sections \ref{SS6} and
\ref{SS7}.  The proof of the Poisson approximation of the loss
network (Theorem \ref{191}) requires some further considerations
presented in the final Section \ref{S62}.

Our work was motivated in part by the posthumous review of Roland
Dobrushin (1996a), where he complained that while ``perturbation
methods are intensively used by mathematical physicists, they are not
so popular as correlation inequalities among the probabilists''.  He
called for a ``systematical exposition oriented to the
mathematicians''.  In this paper we follow his call: we concentrate on
probabilistic issues and exploit probabilistic arguments.

\section{Definitions and results}\label{SS2}

\subsection{The contour model} 
We consider the $d$-dimensional lattice $\Z^d$ and call
\emph{plaquettes} the $(d-1)$-dimensional unit cubes centered at
points of $\Z^d$. We identify each plaquette with its center. A set of
plaquettes is called \emph{surface}. Two plaquettes are
\emph{adjacent} if they share a $(d-2)$-dimensional face. This defines
a notion of connection: a set $\ga$ of plaquettes forms a
\emph{connected} (hyper) surface if for every two plaquettes $x$,$y$ there is
a sequence of pairwise adjacent plaquettes starting at $x$ and ending at $y$.
A surface is \emph{closed} if every $(d-2)$-dimensional face is shared by an
even number of plaquettes in the surface. A \emph{contour} $\ga$ is a
connected and closed family of plaquettes. Two contours $\ga$ and $\th$ are
incompatible if they share some $(d-2)$-dimensional face. In this case we
denote $\ga\not\sim\th$.

Denote by $\G^\Lambda$ the set of all
contours in the volume $\Lambda$.
The subset $\X^\Lambda \subset \N^{\G^\Lambda}$ of \emph{compatible}
configurations is defined as
\begin{equation}
  \X^\Lambda = \{ \eta \in \{0,1\}^{\G^\Lambda}\,;\, \eta(\ga)\,\eta(\theta) =
  0 \mbox{ if } \ga \not\sim\theta\} \label{eq:X}
\end{equation}
that is, a configuration of contours is compatible if it does not contain two
incompatible contours. We denote $\X=\X^{\Z^d}$ and $\G=\G^{\Z^d}$. As usual,
we endowed $\X$ with the product topology.

For each fixed $\beta\in\R^+$, a parameter usually called the inverse
temperature, and for each finite $\Lambda$ define the measure $\mu^\Lambda$ on
$\X^\Lambda$ by
\begin{equation}
  \label{141}
  \mu^\Lambda(\eta) = {\exp\Bigl(-\beta \sum_{\ga:\eta(\ga)=1} |\ga|\Bigr)\over
  Z^\Lambda}
\end{equation}
where $|\ga|$ is the number of plaquettes in $\ga$ and $Z^\Lambda$ is a
renormalization constant making $\mu^\Lambda$ a probability.

\subsection{The loss network}

We introduce a Markov process called (fixed routing) \emph{loss
  network} in the set of compatible contours. This process was
introduced by Erlang in 1917 (see Brockmeyer, Halstrom and Jensen
(1948), p. 139). An account of its properties can be found in Kelly
(1991).  In the traditional interpretation a contour $\gamma$
represents the route taken up by a call. The plaquettes encompasing
$\gamma$ are the circuits held by the call.  For a finite or infinite
set $\Lambda\subset\Z^d$ and $f$ real continuous function on
$\X^\Lambda$, the generator of the process is defined by:
\begin{equation}
\label{204}
A^\Lambda f(\eta) = \sum_{\ga \in \G^\Lambda} \rate\, \one \{\eta+\de_\ga \in
\X^\Lambda\}\, [f(\eta+\de_\ga) - f(\eta)] + \sum_{\ga \in \G^\Lambda}
\eta(\ga)\,[f(\eta-\de_\ga) - f(\eta)]
\end{equation}
where $\de_\ga(\th) = \one\{\th=\ga\}$ and the sum of configurations
is defined pointwisely in $\N^{\G^\Lambda}$:
$(\eta+\xi)(\ga)=\eta(\ga)+\xi(\ga)$.  In words, each contour $\ga$
attempts to appear at rate $e^{-\be |\ga|}$ but it does so only if it
is compatible with all present contours. Present contours disappear at
rate $1$.  Our first result is the following sufficient condition for
the existence of a process with generator~$A^\Lambda$ for any
(infinite) $\Lambda$.
\begin{thm} 
\label{18} If
   \begin{equation} \label{19} \alpha(\beta) = \sup_\ga {1\over
   |\ga|}\,\sum_{\th:\th\not\sim\ga} |\th|\,e^{-\be\vert \th\vert}<
   \infty \end{equation} then for any (infinite) $\Lambda$ the Markov
   process with generator $A^\Lambda$ exists and admits at least one
   invariant measure.
\end{thm}
This theorem is proven in Section \ref{SS6}.  We denote
$\eta^{\Lambda,\zeta}_t$ the corresponding process in $\Lambda$ with
initial configuration $\zeta$.  We omit the volume superindex when
$\Lambda = \Z^d$:
\begin{equation}
  \label{509}
   \eta^\zeta_t = \eta^{\Z^d,\zeta}_t\ \ \ \ A = A^{\Z^d}
\end{equation}

\subsection{Results on the invariant measure} 

We say that $f$ has \emph{support} in $\Upsilon\subset \Z^d$ if $f$
depends only on contours intersecting $\Upsilon$ (not necessarily
contained in $\Upsilon$). Let $|\supp(f)|= \min\{|\Upsilon|: f$ has
support in $\Upsilon\}$. When we write $\supp(f)$ we mean some
$\Upsilon$ such that $|\Upsilon|=|\supp(f)|$ and $f$ has support in
$\Upsilon$. For instance, if $f(\eta) = \eta(\ga)$, $\supp(f)$ may be
set as $\{x\}$ for any $x\in\ga$. The results below work for any such
choice, so one can take the most favorable one in each case. 
 Let
\begin{equation}
 \label{190} \beta^* \hbox{ solution of }\;
\alpha(\beta) = 1,
 \end{equation}
\begin{equation}
  \label{509a}
  \alpha_0(\beta)=\sum_{\ga\ni 0} |\ga|\,e^{-\beta|\ga|}.
\end{equation}
Let $|x|$ be some fixed norm, for instance the one given by the
Manhattan distance ($|x|=\sum_{i=1}^d |x_i|$, the sum of the
coordinate lengths).  Let the corresponding distance between two
subsets of $\Z^d$ be
\begin{equation}
  \label{600}
  d(\Lambda,\Upsilon) = \min \{|x-y|:x\in\Lambda,y\in\Upsilon\}.
\end{equation}

\begin{thm} \label{18a} If $\be>\be^*$ (that is, $\alpha(\be)<1$)  then the
  following statements hold:
  \begin{enumerate}
  \item Uniqueness. For any $\Lambda\subset \Z^d$ there is a unique
    process $\eta_t^\Lambda$ with generator $A^\Lambda$.  The process
    has a unique invariant measure denoted $\mu^\Lambda$.  For
    $\Lambda$ finite this measure is precisely \reff{141}.  For
    $\Lambda=\Z^d$ we denote $\eta_t=\eta_t^{\Z^d}$, $\mu=\mu^{\Z^d}$.
       
  \item Exponential time convergence. For any $\Lambda\subset \Z^d$ and for
    measurable $f$ on $\X^\Lambda$,
    \begin{equation}
      \label{334}
      \sup_{\zeta\in\X^\Lambda} \left\vert \mu^\Lambda f - \E
      f(\eta^{\Lambda,\zeta}_t)\right\vert \;\le\;
      2\,\|f\|_\infty\,|\supp(f)|\,{\alpha_0\over \rho}\, e^{-\rho t}
    \end{equation}
     where $\rho = {1-\alpha\over 2-\alpha}$.
     
   \item Exponential space convergence. Let $\Lambda$ be a (finite or
     infinite) subset of $\Z^d$ and $f$ a measurable function depending on
     contours contained in $\Lambda$.  Then
    \begin{equation}
      \label{130}
      |\mu f - \mu^\Lambda f| \;\le\; 2\,\|f\|_\infty\,\alpha_0\, M_2
      \sum_{x\in\supp(f)}e^{-M_3\, d(\{x\},\Lambda^c)}
    \end{equation}
    where $M_2 = (1-\alpha(\widetilde\beta))^{-1}$ and $M_3 = ( \be -
    \widetilde\be)$, for any $\widetilde\be\in(\beta^*,\beta)$.
    
  \item Exponential mixing. For measurable functions $f$ and $g$ depending
    on contours contained in an arbitrary set $\Lambda\subset\Z^d$:
    \begin{eqnarray}
      \label{101}
      \vert \mu^\Lambda (f g) - \mu^\Lambda f\, \mu^\Lambda g \vert \;\le\;
      2\,\|f\|_\infty\,\|g\|_\infty\,(M_2)^2 \,\sum_{\scriptstyle x\in
        \supp(f),\atop\scriptstyle y\in \supp(g)} |x-y| e^{-M_3 |x-y|}\\[-12pt]
      \nonumber
    \end{eqnarray}
    where $M_2$ ad $M_3$ are the same of \reff{130}.
  \item Central limit theorem. Let $f$ be a measurable function on
    $\X$ with finite support such that $\mu f = 0$ and $\mu(
    |f|^{2+\delta})<\infty$ for some $\delta>0$. Let $\tau_x$ be the
    translation by $x$ and assume $D= \sum_x \mu(f\tau_xf)>0$. Then
    $D<\infty$ and
    \begin{equation}
      \label{90}
      {1\over \sqrt {|\Lambda|}}\, \sum_{x\in \Lambda} \tau_x f
      \sobre{\textstyle\Longrightarrow}{\textstyle\Lambda\to\Z^d} \hbox{\rm
      Normal} (0,D)
    \end{equation}
    where the double arrow means convergence in distribution.
  \end{enumerate}
\end{thm}

This theorem strengthens the results {\bf R1}--{\bf R5} of
Fern\'andez, Ferrari and Garcia (1998). In that paper $\be^*$ was
replaced by a value $\be_M$ defined as the solution of
$\alpha_0(\be_M)=1/(d-1)$.  This value is strictly bigger than
$\be^*$. Item 5 generalizes (the central limit) Theorem 7.4 of
Dobrushin (1996a), where only functions depending on a finite number
of contours are considered. The above theorem will be proven in
Section \ref{SS7}.
   
Finally, we prove a Poisson approximation. Consider the equivalence relation
induced by the translation of contours. Let $\widetilde \G_j$ be a the set
formed by one representative containing the origin from each class of
equivalence of contours with length $j$. For each Borel set $V\subset \R^d$
and $a\in\R$ let
\begin{equation}
  \label{121}
  V\cdot a = \biggl\{x\in \Z^d: \Bigl[{x-1/2\over a},{x+1/2\over
    a}\Bigr]^d\subset V\biggr\}.
\end{equation}
Fix a contour length $j>0$. For each $\ga\in\widetilde \G_j $, let
$M_{\ga,\beta}$ be defined by
\begin{equation}
\label{920}
M_{\ga,\beta}(V) = \sum_{x\in V\cdot e^{\be |\ga|/d}} \eta^\beta(\tau_x\ga) 
\end{equation}
where $\eta^\be$ is distributed according to the invariant measure for the
given $\be$. Let $(M_{\ga,\infty},\, \ga\in\widetilde \G_j )$ be a family of
independent unit Poisson processes in $\R^d$.

\begin{thm}
      \label{191}
      For each contour length $j$, it is possible to jointly construct the
      $d$-dimensional processes
      $\{M_{\ga,\beta}:\be^*<\be\le\infty,\,\ga\in\widetilde \G_j\}$ in such a
      way that for all regions $V$ equal to a product of intervals,
\begin{equation}
  \label{125}
  \P\Bigl(M_{\ga,\beta}(V)\neq M_{\ga,\infty}(V)\Bigr) \; \leq \;
  c(|\ga|,V)\Bigl(\alpha(\beta)+e^{-\be|\ga|/d}\Bigr)\; \sim \;
  \exp\Bigl(-\be\, \min\{2d, |\ga|/d\}\Bigr).
\end{equation}
where $c(|\ga|,V)$ is a computable constant.  As a consequence $
\{M_{\ga,\beta}:\ga\in\widetilde \G_j\}$ converges in distribution to a family
consisting of $|\widetilde \G_j|$ independent Poisson processes with mean
$e^{-\be j}$.
 \end{thm} 
 
 Theorem \ref{191} is proven in Section \ref{S62}.

\section{Graphical representation of loss networks}\label{S3}

We construct the loss network as a function of stationary marked
Poisson processes (a la Harris) each of which indicates the attempted
birth-times of a contour. A lifetime is associated to each attempted
birth. The triple (contour , attempted birth , lifetime) is called a
\emph{cylinder}. The loss network is constructed by erasing cylinders
which at birth violate the exclusion condition. The crucial point in
this construction is the association of the lifetime to the
birth-time. This allows us to study the process backwards in time by
studying a Markovian oriented percolation process of cylinders. In
contrast, the standard construction uses independent Poisson processes
for the birth-times and death-times respectively. In this case the
backwards construction looks hard.

\subsection {Marked Poisson processes} \label{4.1}

To each contour $\ga \in \G$ we associate an independent (of everything)
marked Poisson process $N_{\ga}$ with rate $e^{-\be |\ga|}$. We call
$T_{k}(\ga)$, $\ga\in\G$, the ordered time-events of $N_\ga$ with the
convention that $T_{0}(\ga)<0<T_{1}(\ga)$.  For each occurrence time
$T_{i}(\ga)$ of the process $N_{\ga}$ we choose an independent mark
$S_{i}(\ga)$ exponentially distributed with mean 1. At the Poisson time-event
$T_i(\ga)$ a contour $\ga$ appears and it lasts $S_i(\ga)$ time units.
  
The random family $\C=\bigl\{\{(\ga, T_{i}(\ga),S_{i}(\ga)): i \in \Z
  \}:\ga\in \G\bigr\}$ consists of independent marked Poisson processes.
A marked point  $(\ga,T_k(\ga),S_k(\ga))\in\C$ is identified with $\ga \times
[T_k(\ga), T_k(\ga)+S_k(\ga)]$, the \emph{cylinder} with \emph{basis} $\ga$,
\emph{birth-time} $T_k(\ga)$ and  \emph{lifetime}  $S_k(\ga)$. The \emph{life}
of the cylinder is the time interval
  $[T_k(\ga),T_k(\ga)+S_k(\ga)]$. For a generic cylinder
$C=(\ga,t,s)$, we use the notation
\begin{equation}
  \label{144}
  \basis(C) = \ga,\ \ \birth(C)= t,\ \ \death(C)= t+s,\ \ \life(C) = [t,t+s].
\end{equation}
We define incompatibility between cylinders $C$ and $C'$ by
\begin{equation} 
\label{502}
C'\not\sim C \;\;\hbox{ if and only if }\;\;
\basis(C)\not\sim\basis(C')\; \hbox{ and
  }\;\life(C)\cap\life(C')\neq\emptyset, 
\end{equation}
otherwise $C'\sim C$ (compatible). We say that two sets of cylinders
$\A$ and $\A'$ are \emph{incompatible} if there is 
a cylinder in $\A$ incompatible with a cylinder in $\A'$:
\begin{equation}
  \label{230}
  \A \not\sim \A'\;\hbox{ if and only if } \;C\not\sim C' \hbox{ for some }
  C\in \A \hbox{ and } C'\in \A'.
\end{equation}
  
Let $\bS:=(S^0_i(\th):\th\in\G,i\ge1)$ be a (countable) family of iid
exponential times of mean $1$ independent of $\C$. These are the
lifetimes that, when necessary, will be associated to the contours of
the initial configuration. Indeed, we identify $\bS$ with the set of
cylinders $\{(\th,0,S^0_i(\th)):\th\in\G, i\ge 1\}$. For $\xi\in
\N^{\G}$, let
\begin{equation}
  \label{c00}
  \bS(\xi)= \bigcup_{\th\in\G}\bigcup_{i=1}^{\xi(\th)}\{(\th,0,S^0_i(\th))\}
\end{equation}
the family of cylinders associated to the initial configuration $\xi$,
all with birth-time $0$. Notice that $\xi$ may have more than one
cylinder with the same basis.
For $s<t$ define
\begin{equation}
  \label{cst}
  \C{[s,t]} := \{C\in \C: \birth(C)\in [s,t]\}.
\end{equation}
the set of cylinders born in the interval $[s,t]$.
 
{\bf Remark.} In this paper we will work with the probability space
given by the product of the spaces generated by $\C$ and $\bS$. We
call it $(\Omega,{\cal F},\P)$. We write $\E$ for the respective
expectation. In section \ref{S30} we use the direct product of this
space by itself, while in Section \ref{SSS} we need to consider
countable products of this space. We use the same notation
$\P$ and $\E$ for the corresponding probability and expectation in
these enlarged spaces.

\subsection{The free network}\label{free}

For $\xi\in\N^\G$ define 
\begin{equation}
  \label{201}
  \xi^\xi_t(\ga) = \sum_{C\,\in\,\C[0,t]\cup\bS(\xi)} \one\{\basis(C)=\ga,
  \life(C)\ni t\}. 
\end{equation}
The above process, called the \emph{free network}, is a product of
independent birth-and-death processes on $\N^{\G}$ with initial
configuration $\xi$ whose generator is given by
\begin{equation}
  \label{144b}
A^{0} f(\xi) = \sum_{\ga \in \G} \rate [f(\xi+\de_\ga) -
f(\xi)] + \sum_{\ga \in \G} \xi(\ga)[f(\xi-\de_\ga) - f(\xi)].
\end{equation}
The invariant (and reversible) measure for this process is the product measure
$\mu^{0}$ on $\N^\G$ with Poisson marginals
\begin{equation}
\label{eq:muga}
\mu^{0}\bigl\{\xi(\ga)=k\bigr\} \;=\;  
\frac{(\rate)^{k}}{k!} \,\exp\left(\rate\right)\;.
\end{equation}

In terms of loss networks, $\xi^\xi_t$ is the process for which all
the calls are accepted, that is $\xi^\xi_t(\ga)$ is the number of
calls on route $\ga$ at time $t$ when there is no restriction on the
number of calls a circuit can accept and the initial configuration of
calls is $\xi$.

\subsection{Finite-volume construction of a loss network}\label{3.3}

In the construction of a loss network in a finite volume $\Lambda$
with an initial condition $\zeta\in\X^\Lambda$ we use only the finite
set of Poisson processes $(N_\ga:\ga\subset\Lambda)$ and the finite
family of initial lifetimes $(S^0_1(\th): \th\subset\Lambda)$.  Let
$\C^\Lambda = \{C\in\C:\basis(C) \subset \Lambda\}$ and $\bS(\zeta)$,
defined as in \reff{c00}, be such that all its cylinders are mutually
compatible.
We realize the dynamics $\eta^{\Lambda,\zeta}_t$ as a
(deterministic) function of $\C^\Lambda$ and $\bS(\zeta)$.

We construct inductively $\Ke^\Lambda[0,t]$, the set of \emph{kept}
cylinders at time $t$. The complementary set corresponds to
\emph{erased} cylinders. At time zero we include all cylinders of
$\bS(\zeta)$ in $\Ke^\Lambda[0,t]$.  Then, we move forward in time
and consider the first Poisson mark: The corresponding cylinder is
erased if it is incompatible with any of the cylinders already in
$\Ke^\Lambda[0,t]$, otherwise it is kept.  This procedure is
successively performed mark by mark until all cylinders born before
$t$ are considered.  Define $\eta^{\Lambda,\zeta}_t\in\X^\Lambda$ as
\begin{eqnarray}
  \label{12}
 \eta^{\Lambda,\zeta}_t(\gamma)  &=& 
\sum_{C\in\Ke^\Lambda[0,t]} \one\{\basis(C)=\ga,
  \life(C)\ni t\}\\ 
  &=& \one\Bigl\{ \ga \in
  \bigl\{\basis(C): C\in\Ke^\Lambda[0,t],\, \life(C)\ni t\bigr\}\Bigr\},
\end{eqnarray}
that is, $\eta^{\Lambda,\zeta}_t$ signals all contours which are basis
of a kept cylinder that is alive at time $t$.  We show in Section
\ref{S34} that $\eta^{\Lambda,\zeta}_t$ has generator $A^\Lambda$
defined as in \reff{204} restricting the sums to the set of contours
contained in $\Lambda$. It is immediate that $\mu^\Lambda$ defined in
\reff{141} is reversible for this process. Since we are dealing with
an irreducible Markov process in a finite state space,
$\eta^{\Lambda,\zeta}_t$ converges in distribution to $\mu^\Lambda$
for any initial configuration $\zeta$. This in particular implies that
$\mu^\Lambda$ is the unique invariant measure for this process. Later
in the paper we determine the speed of convergence.

Using the same $\C$ and $\bS$ in the construction of $\eta_t^{\Lambda,\zeta}$
and $\xi^\xi_t$, we have that if $ \zeta(\ga)\;\le\;
  \xi(\ga)$ for all $\ga\subset\Lambda$, then
\begin{equation}
\label{621a}
  \eta_t^{\Lambda,\zeta}(\ga)\;\le\;
  \xi^\xi_t(\ga), \;\;\hbox{for all $\ga\subset\Lambda$,}
\end{equation}
because in the free network $\xi^\xi_t$ all cylinders are kept. 

Since $\Lambda$ is finite, there exists a sequence of random times
$t_i=t_i(\C^\Lambda)$ with $t_i\to\pm\infty$ as $i\to\pm\infty$ such
that $\xi_{t_i}(\ga) = 0$ for all $\ga\in\Lambda$. We can, in
particular, consider $t_i$ as the entrance times of $\xi_t$ in the set
$\{\xi:\xi(\gamma)=0,\;$ for all $\gamma\in\Lambda\}$. Since this
process has a unique invariant measure which gives positive
probability to this set, $(t_i)$ is a stationary renewal process with
inter-renewal time with finite mean. We extend the construction of a
set of kept cylinders to $t\in\R$ forgetting the set $\bS(\zeta)$, by
doing the above procedure in each time interval $[t_i,t_{i+1}]$ with
the cylinders of $\C[t_i,t_{i+1}]$. This can be done because no
cylinder intersects $\{t_i:i\in\Z\}$. Let us denote $\K^\Lambda$ the
resulting set of kept cylinders and $\eta^\Lambda_t$ its projection in
the sense of \reff{12}. By construction $\K^\Lambda$ has a
time-translation invariant distribution. The process $\eta^\Lambda_t$
has generator $A^\Lambda$ and distribution independent of $t$, hence
given by $\mu^\Lambda$. This implies that for any
$f:\{0,1\}^{\G^\Lambda}\to \R$ and any $t\in\R$,
\begin{equation}
  \label{506}
  \mu^\Lambda f = \E f(\eta^\Lambda_t).
\end{equation}
Since $\eta^\Lambda_t(\ga) \le \xi_t(\ga)$ for all $\gamma\in\Lambda$
and $\xi_t(\ga)$ has Poisson distribution with mean $e^{-\be|\ga|}$,
we have, taking $f(\eta) = \eta(\ga)$ in \reff{506}, that
\begin{equation}
\label{621}
\mu^\Lambda\{\eta:\eta(\ga)= 1\}\;\le\;
  e^{-\be|\ga|}.
\end{equation}

\subsection{Infinite-volume construction. 
Backwards oriented percolation}\label{S34}

If we try to perform an analogous construction in infinite volume we are
confronted with the problem that there is no first mark. To overcome this
we follow the original approach of Harris (1972) [see also Durrett (1995)] and
introduce the notion of percolation. The goal is to partition the set of
cylinders in finite subsets to which the previous mark-by-mark construction
can be applied.

We come back to the infinite-volume construction of Section \ref{4.1}. For an
arbitrary space-time point $(x,t)$ define the set of cylinders
containing the point $(x,t)$ by
\begin{eqnarray}
  {\bf A}^{x,t}_1
 &=& \{C\in\C\, ;\, \basis(C)\ni x,\ \life(C)\ni t\}\,.  \label{140}
\end{eqnarray}

For any cylinder $C$ define the set of {\it ancestors} of $C$
as the set of cylinders born before $C$ that are incompatible
with $C$:  
\begin{eqnarray}
  {\bf A}^{C}_1 &=& \{C'\in\C\, ;\, C'\not\sim C,
  \birth(C')<\birth(C)
\}
  \label{140a} \\
&=& \bigcup_{x\in\basis(C)} \A^{x,\birth(C)}\,.
\end{eqnarray}
The definition of ancestor of $C$ does not depend on the lifetime
of $C$.  Recursively for $n\ge 2$, the $n$th generation of ancestors
of $(x,t)$ is defined as
\begin{equation}
  {\A}_{n}^{x,t} = \{ C'':C''\in{\A}^{C'}_1\hbox{ for some } C' \in
  {\A}_{n-1}^{x,t} \},  \label{eq:an}
\end{equation}
and for a given cylinder $C$,
\begin{equation}
{\A}_n^C = \{ C'':C''\in{\bf A}^{C'}_1\hbox{ for some } C' \in {\bf
      A}_{n-1}^C \}.\label{eq:ances}
\end{equation}

We say that there is {\it backward oriented percolation} in ${\C}$ if there
exists a space-time point $(x,t)$ such that ${\bf A}_{n}^{x,t} \neq \emptyset$
for all $n$, that is, there exists a point with infinitely many generations of
ancestors.  Let the \emph{clan} of the space-time point $(x,t)$ be the union
of its ancestors:
\begin{equation}
  \label{206}
  \A^{x,t}=\bigcup_{n\ge 1}{\A}_{n}^{x,t}.
\end{equation}

In the next theorem we give a sufficient condition for the existence of the
infinite-volume process in any finite time interval in terms of backwards
percolation.

\begin{thm}
  \label{20} If with probability one $\A^{x,t}\cap\C{[0,t]}$ is finite
  for any $x\in\Z^d$ and $t\ge 0$, then for any (possibly infinite)
  $\Lambda\subset\Z^d$, the process with generator $A^\Lambda$ is well
  defined for any initial configuration $\zeta\in\X^\Lambda$and has
  at least one invariant measure $\mu^\Lambda$.
\end{thm}

\noindent{\bf Proof.}
We construct the process for $\Lambda = \Z^d$; the construction for
other regions $\Lambda$ is analogous. The initial configuration is
denoted $\zeta\in\X$ and the initial cylinders are given by
$\bS(\zeta)$, defined in \reff{c00}.  Note that that all the cylinders
of $\bS(\zeta)$ are mutually compatible. We then partition
$\bS(\zeta)\cup\C[0,t]$ into a set of \emph{kept} cylinders, denoted
by $\K$, and a set of \emph{erased} cylinders denoted by $\D$.

The construction is as follows.  First, all cylinders in $\bS(\zeta)$
are kept.  Second, for each $x\in\Lambda$ the percolation clan of
$(x,t)$ in $[0,t]$, $\A^{x,t}\cap \C[0,t]$, is partitioned in kept and
deleted cylinders as in the finite-volume case.  To do so we order the
cylinders of $\A^{x,t}\cap \C[0,t]$ by birth-time.  This can be done
because by hypothesis $\A^{x,t}\cap \C[0,t]$ has a finite number of
cylinders.  Then we successively classify each cylinder as kept if it
is compatible with all cylinders already classified as kept (including
those in $\bS(\zeta)$); if not, we classify it as erased. We denote the
resulting sets $\Ke^{x,t}[0,t]$ and $\D^{x,t}_{\zeta}[0,t]$,
respectively.

Denoting
\begin{equation}
  \label{k0t}
  \Ke[0,t]\;:=\;\bigcup_{x\in\Lambda}\Ke^{x,t}[0,t]\;,\;\;\;\;\;
  \D_{\zeta}[0,t]\;
:=\;\bigcup_{x\in\Lambda}\D_{\zeta}^{x,t}[0,t]\;,
\end{equation}
we have that
\begin{equation}
  \label{eq:kdc}
  \Ke[0,t]\, \dot\cup\, \D_{\zeta}[0,t] \;=\; \C[0,t]\,\cup\,\bS(\zeta)\;.
\end{equation}
Indeed, the classification of any given cylinder $C\in \C[0,t]$
depends only on (a) its ancestors in $[0,t]$, $\A^C\cap\C[0,t]$, and
(b) on the finite subset of $\bS(\zeta)$ of cylinders which are
incompatible with some of the ancestors of $C$ in $[0,t]$. Therefore
there is no inconsistency: $\Ke^{x,t}[0,t]\,\cap \,
\D_{\zeta}^{y,t}[0,t]\;=\;\emptyset$ for all $x\neq y$ and
$\Ke^{x,t}[0,t]\,\subset\,\Ke^{x',t'}[0,t']$ for $t<t'$, if $(x,t)\in
C$ for some $C\in \Ke^{x',t'}[0,t']$.

The process is now defined as in \reff{12} by
\begin{equation}
  \label{152} \eta^{\zeta}_t(\ga) =\one \Bigl\{\ga\in\bigl\{\basis(C): C\in
  \Ke[0,t],\, \life(C)\ni t\bigr\}\Bigr\}.
\end{equation}
The reader can check that for finite $\Lambda$ the above construction
is equivalent to that of Section \ref{3.3}.  Applied to the
set of cylinders of $\C^\Lambda[0,t]$ it yields the set
$\Ke^\Lambda[0,t]$ defined in the paragraph preceding formula
\reff{12}. 

To show that $\eta^{\zeta}_t$ has generator $A$, denote $\eta_t =
\eta^\zeta_t$ and $\K= \K_\zeta$ and write
\begin{eqnarray}
  \label{153}
  \lefteqn{ [f(\eta_{t+h}) -  f(\eta_t)] }\nonumber\\
  &=& \sum_{C\in \K[0,t+h]} \one\{\birth(C)\in [t,t+h]\}
  [f(\eta_t+\de_{\basis(C)})-f(\eta_t)]\nonumber\\
  &&{}+\sum_{C\in \K[0,t]} \one\{\life(C)\ni t,\, \life(C)\not\ni
  t+h\}[f(\eta_t-\de_{\basis(C)})-f(\eta_t)] \nonumber\\
  &&{}+ \{\hbox{other things}\},
\end{eqnarray}
where $\{${other things}$\}$ refer to events with more
than one Poisson mark in the time interval $[t,t+h]$ for the
contours in the (finite) support of $f$. Since the total rate of the
Poisson marks in this set is finite, the event $\{$other things$\}$ has
a probability of order $h^2$. Now, denoting
\begin{eqnarray}
\label{214}
N_\ga(t,s) &=& \#\{k:T_k(\ga)\in (t,s)\}
\end{eqnarray}
we have
\begin{eqnarray}
  \label{154} \lefteqn{ \sum_{C\in \C} \one\{\birth(C)\in
  [t,t+h]\}\,\one\{C\in
  \K[0,t+h]\}\,[f(\eta_t+\de_{\basis(C)})-f(\eta_t)]}\nonumber\\ &
  =&\sum_{C\in \C} \one \{N_{\basis(C)}[t,t+h]=1\}\,
  \one\{\basis(C)\sim \basis(C'), \forall C'\in\K[0,t]: \life(C')\ni
  t\}\nonumber\\ &&\qquad\qquad\qquad\qquad\qquad\qquad\qquad
  \qquad\qquad\qquad\qquad \qquad\qquad
  \times[f(\eta_t+\de_{\basis(C)})-f(\eta_t)]\nonumber\\ & =&
  \sum_\ga\one \{N_{\ga}[t,t+h]=1\}\,\one\{\eta_t+\de_{\ga}\in\X\}
\,[f(\eta_t+\de_{\ga})-f(\eta_t)]   
\end{eqnarray}

To compute the second term of \reff{153}, observe that $\life(C)$ is
independent of $\birth(C)$ and both the event $\{C\in \K[0,t]\}$ and
$\eta_t$ are ${\cal F}_t$-measurable. Here ${\cal F}_t$ is the
$\sigma$-algebra generated by the births and deaths occurred before
$t$.  Hence
\begin{eqnarray}
\P\Bigl(\life(C)\ni
  t,\, \life(C)\not\ni
  t+h \Given{\cal F}_t\Bigr)
  &=& \P\Bigl(\life(C)\not\ni t+h \Given \life(C)\ni t\Bigr) 
 \label{1555}
\end{eqnarray}
and
\begin{eqnarray}
  \label{155}
 \lefteqn{
 \E \Bigl[ \sum_{C} \one\{C\in\K[0,t]\}\,\one\{\life(C)\ni
  t,\, \life(C)\not\ni
  t+h\}\,[f(\eta_t-\de_{\basis(C)})-f(\eta_t)]\Bigr]
        }
\nonumber \\
 & =&
  \E \Bigl[  \sum_C \P\Bigl(\life(C)\not\ni
  t+h\Given \life(C)\ni t\Bigr)\,\one\{C\in\K[0,t],\;\life(C)\ni
  t\}\,\nonumber \\
&&\qquad\qquad\qquad\qquad\qquad\qquad\qquad\qquad\qquad\qquad\times
  [f(\eta_t -\de_{\basis(C)})-f(\eta_t)]\Bigr].
\end{eqnarray}
Since $\life(C)$ is exponentially distributed with mean 1,
\begin{eqnarray}
\P\Bigl(\life(C)\not\ni t+h \Given \life(C)\ni t\Bigr) 
&=&  h + o(h)\,. \label{1554}
\end{eqnarray}
Taking the expectation of \reff{153} and substituting
\reff{154}--\reff{1554} we get
\begin{eqnarray}
  \label{156} \lefteqn{\E [f(\eta_{t+h}) - f(\eta_t)] }\nonumber\\ &
  = & \sum_\ga h \,\rate \,\E\Bigl(\one\{\eta_t+\de_{\ga}\in\X\}\,
    [f(\eta_t+\de_{\ga})-f(\eta_t)]\Bigr) + o(h)\nonumber \\ & & + \sum_\ga
  h\,\E \Bigl(\eta_t(\ga) \,[f(\eta_t-\de_{\ga})-f(\eta_t)] \Bigr)+o(h)
\end{eqnarray}
which dividing by $h$ and taking limit gives 
\begin{equation}
  \label{dge}
  {d\E f(\eta^\zeta_t)\over dt} \;=\; A \E f(\eta^\zeta_t).
\end{equation}

The existence of an invariant measure follows by compactness as our
process is defined in the compact space $\X$. See  Chapter 1 of
Liggett (1985). \   \square

We show in the next theorem that under stronger hypothesis the
process can be constructed for times in the whole real line. Since the
construction is time-translation invariant, the distribution of
$\eta_t$ will be invariant.

\begin{thm}
  \label{20a}
  If with probability one there is no backwards oriented percolation in $\C$,
  then the process with generator $A$ can be constructed in $(-\infty,\infty)$
  in such a way that the marginal distribution of $\eta_t$ is invariant.
\end{thm}

\proof The lack of percolation allows us to construct a set $\K\subset
\C$ as $\K_\zeta[0,t]$ was constructed from $\C[0,t]\cup\bS(\zeta)$ in
the proof of the previous theorem. We just proceed clan by clan and
simply ignore the cylinders of $\bS$.  Note that $\K$ is both space
and time translation-invariant by construction.  Analogously to the
previous theorem we define $\eta_t$ as the section of $\K$ at time
$t$:
\begin{equation}
  \label{210}
  \eta_t(\ga)  =\one \Bigl\{\ga\in\bigl\{\basis(C): C\in \K,\,
  \life(C)\ni t
\bigr\}\Bigr\}.
\end{equation}
By construction, the distribution of $\eta_t$ does not depend on $t$, hence
its distribution is an invariant measure for the process.  \square
\bigskip

Let us denote $\mu$ the distribution of $\eta_t$, in anticipation of
the fact that this is precisely the measure of Theorem \ref{18a}.
\bigskip

As in the finite case, 
\begin{equation}
  \label{a623}
  \eta_t(\ga)\;\le\; \xi_t(\ga)
\end{equation}
for all $\ga\in \G$.  This implies that the distribution $\mu$
inherits property \reff{621}:
\begin{equation}
  \label{623}
  \mu\bigl(\eta_t^\Lambda(\ga)= 1\bigr) \;=\; \E\eta_t(\ga)\;\le\;
  \E\xi_t(\ga) \;=\; e^{-\be|\ga|}.
\end{equation}

Let 
\begin{equation}
  \label{213}
  \A(\Upsilon) = \cup_{x\in\Upsilon}\A^{x,0}
\end{equation}
be the \emph{clan} of $\Upsilon\subset\Z^d$ (at time $0$).
\medskip

\noindent{\bf Remarks.}\
{\bf 1.}\ It follows from \reff{210} that for any $t\in \R$ and continuous $f$,
\begin{equation}
  \label{505}
  \mu f = \E f(\eta_t).
\end{equation}

{\bf 2.}\ The presence/absence of a contour $\ga$ at time $t$ depends
only on the clan of ancestors of $(x,t)$ for any $x\in\ga$ through a
certain function.  More generally, for each $f$ there exists a
function $\Phi$ such that
\begin{equation}
  \label{77}
  f(\eta_t) = \Phi\Bigl(\cup_{x\in \supp(f)}\A^{x,t}\Bigr).
\end{equation}
For instance, $\eta_t(\ga)=1$ if and only if $\A_1^{x,t}$ contains a cylinder
in $\K$ with basis $\ga$ whose life contains $t$. This depends only on the set
$\A^{x,t}$. In particular, with the notation \reff{213},
\begin{equation}
  \label{77a}
  f(\eta_0) = \Phi\Bigl(\A(\supp(f))\Bigr).
\end{equation}
Analogous statements are true for the process starting with a fixed
configuration at time zero:
\begin{equation}
  \label{77k}
  f(\eta^\zeta_t) = \Phi\Bigl(\A_\zeta( \supp(f),[0,t])\Bigr)
\end{equation}
where
\begin{equation}
  \label{az0}
 \A_\zeta(\Upsilon,[0,t])\;=\; \Bigl[(\cup_{x\in\Upsilon}
  \A^{x,t})\, \cap\,   \C[0,t]\Bigr]\,\cup\,
\Bigl\{C\in\bS(\zeta): \{C\}\,\not\sim\,
   (\cup_{x\in\Upsilon} \A^{x,t})\, \cap\,
  \C[0,t] \Bigr\}
\end{equation}
is the set of cylinders in $\C[0,t]\cup\bS(\zeta)$ which determines 
the value of $f(\eta^\zeta_t)$ when $f$ has support~$\Upsilon$.

\section{Percolation, space-time convergence and mixing}\label{SS4}

In this section we exploit the relation between the loss-network
process and the absence of percolation to prove a more precise version
of Theorem \ref{18a}. In the proof of the mixing properties we shall
need a continuous-time construction of the backwards percolation clan.

\subsection{The key theorem}\label{s130}

The precise statement of next theorem requires the notion of
non-oriented percolation in a time interval. For any time interval
$(s,t)$ and any space-time point $(x,t')$ define
\begin{equation}
  \label{221}
  {\G}_{0}^{x,t'}[s,t]= \Bigl\{ C\in\C[s,t]: \basis(C)\ni x, 
\life(C)\ni t' \Bigr\}
\end{equation} 
and
\begin{equation}
  {\G}_{n}^{x,t'}[s,t] = \Bigl\{ C\in\C[s,t]: \basis(C)\not\sim \basis(C')
  ,\,\hbox{ for some } C'\in  {\G}_{n-1}^{x,t'}[s,t] \Bigr\}.
  \label{220}
\end{equation}
Notice that in the definition of $\G_n$ there is no exigency that the birth
time of $C'$ be previous to the birth-time of $C$ or that the lifetimes
intersect.  Let
\begin{equation}
  \label{222}
  \G^{x,t'}[s,t] = \bigcup_{k\ge 0} {\G}_{k}^{x,t'}[s,t].
\end{equation}
We say that there is no (non-oriented) percolation in $[s,t]$ if for
any space-time point $(x,t')$, $\G^{x,t'}[s,t]$ contains a finite
number of cylinders. We will show later that the condition
$\alpha<\infty$ is sufficient for the existence of an $h$ such that
the probability that there is no non-oriented percolation in $[0,h]$
is one.

In addition we need the following definitions.  
\begin{itemize}
\item
The \emph{time-length}
and the \emph{space-width} of the family of cylinders $\A^{x,t}$ are
respectively
\begin{eqnarray}
  \label{229}
  \tl(\A^{x,t}) &=& t- \sup\{s: \life(C)\ni s, \hbox{ for some } C\in
  \A^{x,t}\},\\
 \sw(\A^{x,t}) &=& \Bigl|\bigcup_{C\in \A^{x,t}} \basis(C)\Bigr|.\label{232}
\end{eqnarray}
In words, the space-witdth is the number of sites occuppied by the
projection of the bases of the cylinders in the family. The
time-lenght is the lenght of the time interval between $t$ and the
first birth in the family of ancestors of $(x,t)$.

\item $\A^\Lambda(\supp(f))$ is the set of ancestors of $\supp(f)$
constructed from $\C^\Lambda$ as $\A(\supp(f))$ was constructed from
$\C$. [Notice that this is \emph{not} the same as
$\A(\supp(f))\cap\C^\Lambda$].

\end{itemize}

In item 4 of next theorem we enlarge our probability space to the
direct product of our working space with itself: $(\Omega,{\cal
F},\P)\times (\Omega,{\cal F},\P)$. As warned before, we continue to
use $\P$ and $\E$ for the probability and expectation of this space.

\begin{thm}
  \label{225}
  Assume that there is no backwards oriented percolation with probability
  one. 
  Then,
\begin{enumerate}
\item Uniqueness. The measure $\mu$ is the unique invariant measure for the
    process $\eta_t$.

  \item Time convergence. For any function $f$ with finite support,
    \begin{equation}
      \label{754}
   \lim_{t\to\infty}\sup_{\eta\in\X} \left|\E f(\eta^{\zeta}_t)-
    \mu f\right|=0.
    \end{equation}
Furthermore,
 \begin{eqnarray}
   \lefteqn{\sup_{\eta\in\X} \left| \mu f - \E f(\eta^{\zeta}_t)\right|}\nonumber\\
   &\le& 2\,\|f\|_\infty\,\P\Bigl(\cup_{x\in\supp(f)}\{\A^{x,t} \not\sim
   \bS(\zeta) \hbox{ or } \tl(\A^{x,t})>t\}\Bigr)\label{751}\\
   &\le& 2\,\|f\|_\infty\,\sum_{x\in\supp(f)}\Bigl[\P\Bigl( \tl(\A^{x,0})>bt \Bigr) +
   e^{-(1-b)t}\, \E\Bigl(\sw(\A^{x,0})\Bigr)\Bigr] \label{335}
   \end{eqnarray}
for any $b\in(0,1)$.

\item Space convergence. As $\Lambda\to\Z^d$, $\mu^\Lambda$ converges weakly
  to $\mu$. More precisely, if $f$ is a function depending on contours
  contained in a finite set $\Lambda$, then
\begin{eqnarray}
  {|\mu f- \mu^\Lambda f|}\;\le
  \;2\,\|f\|_\infty\,\P\Bigl(\A(\supp(f)) \neq
  \A^\Lambda(\supp(f))\Bigr)\, .\label{73a}
\end{eqnarray}

\item Mixing. If in addition there exists a value $h$ such that there is no
  (non-oriented) percolation in $(0,h)$ with probability one, then for $f$ and
  $g$ with finite support, 
  \begin{equation}
    \label{73d}
    \lim_{|x|\to\infty} \vert \mu (f \tau_x g) - \mu f\, \mu g \vert\; =\; 0
  \end{equation}
where $\tau_x$ is translation by $x$. More precisely,
\begin{equation}
 \label{73b}
{|\mu (fg)- \mu f\, \mu g|}
 \;\le\; 2\,\|f\|_\infty\|g\|_\infty\, \P\Bigl(\A(\supp(f)) \, \not\sim\,
  \widehat\A(\supp(g))\Bigr)
\end{equation}
where $\widehat \A(\supp(g))$ has the same distribution as $\A(\supp(g))$
but is independent of\break $\A(\supp(f))$.
\end{enumerate}
\end{thm}

Existence of $\mu$ has been proven in Theorem \ref{20a}.  In the rest of the
section we prove the other properties.

\subsection{Time convergence and uniqueness} \label{tcu}

We use the same Poisson marks to construct simultaneously the
stationary process $\eta_t$ and a process starting at time zero with
an arbitrary initial configuration $\zeta$. The second process is
denoted $\eta^{\zeta}_t$ (as before). This is what in the literature
is known as \emph{coupling}. By construction [cf.\/ \reff{k0t} and
\reff{152}], the process $\eta^{\zeta}_t$ ignores the cylinders in
$\C$ with birth-times less than $0$ but takes into account the set of
cylinders with basis given by the contours of the initial
configuration $\zeta$ and birth-time zero,
$\bS(\zeta)=\{(\th,0,S^0_1(\th))\in \bS:\eta(\th)=1\}$. Recall that
the times $S^0_k(\th)$ are exponentially distributed with mean 1 and
independent of everything.

By \reff{505} and \reff{506},
\begin{equation}
\label{34b}
\sup_{\zeta\in\X}\,\Bigl|\E f(\eta^\zeta_t)- \mu f\Bigr|\;=\;
\sup_{\zeta\in\X}\, \left\vert\E\Bigl(f(\eta^\zeta_t) -
  f(\eta_t)\Bigr)\right\vert .  
\end{equation}
Since we are using $\C$ to construct $\eta_t$ and
$\C[0,t]\cup \bS(\zeta)$ to construct $\eta^\zeta_t$, it follows from
\reff{77} and \reff{77k} that
\begin{eqnarray}
  \nonumber
  {\Bigl| f(\eta^\zeta_t) - f(\eta_t) \Bigr| }&=&\Bigl|\Phi\Bigl(\A_\zeta
  ( \supp(f),[0,t])\Bigr)- 
 \Phi\Bigl(\cup_{x\in\supp(f)} \A^{x,t}\Bigr)\Bigr|\\
&\le&
2\,\|f\|_\infty\,\one\left\{\Bigl(\bS(\zeta)\not\sim\cup_{x\in\supp(f)}
  \A^{x,t} \Bigr)
    \hbox{ or } \tl(\A^{x,t})>t\right\}\label{227}
\end{eqnarray}
To see this notice that $\{\A^{x,t} \subset
\C[0,t]\}=\{\tl(\A^{x,t})<t\}$, and that $\cup_{x\in\supp(f)}
\A^{x,t}\sim \bS(\zeta)$ and $\cup_{x\in\supp(f)} \A^{x,t} \subset
\C[0,t]$ if and only if $\cup_{x\in\supp(f)} \A^{x,t}\,= \,\A_\zeta (
\supp(f),[0,t])$ .  Display \reff{227} shows \reff{751}.

To prove the weak convergence \reff{754} we fix $b\in[0,1]$ and bound the
indicator function in the right hand side of \reff{227} by
\begin{equation}
  \label{228}
  \one\Bigl\{ \tl(\A^{x,t})>bt\Bigr\} + 
\one \Bigl\{\tl(\A^{x,t})<bt , \A^{x,t} \not\sim
  \bS(\zeta)\Bigr\}.
\end{equation}
The expected value of the first term in \reff{228} goes to zero
because $\A^{x,t}$ has a finite number of cylinders with probability
one. The second term in \reff{228} is bounded above by
\begin{equation}
  \label{2288}
  \one \Bigl\{\max \{S^0_1(\ga): \zeta(\ga)=1\hbox{ and }
\ga\not\sim\basis(C)
  \hbox{ for some }C \in \A^{x,t}\}>(1-b)t\Bigr\}.
\end{equation}
Since $\bS$ and $\A^{x,t}$ are independent and $S^0_i$ are iid
exponentially distributed random variables of mean $1$,
\begin{eqnarray}
  \lefteqn{\hspace{-2cm}\E\Bigl(\max \{S^0_1(\ga): \zeta(\ga)=1\hbox{ and }
  \ga\not\sim\basis(C)
  \hbox{ for some }C \in \A^{x,t}\}>(1-b)t\,\Big\vert\,
  \A^{x,t}\Bigr)}\nonumber\\
&&= 1-(1-e^{-(1-b)t})^{|\{\ga: \zeta(\ga)=1\hbox{ and }
\ga\not\sim\basis(C)
  \hbox{ for some }C \in \A^{x,t}\}|}\;.  \label{2289}
\end{eqnarray}
Since $\zeta$ is a configuration of compatible contours, it contains
at most one contour per site, i.e.\ $|\{\ga\ni x:\zeta(\ga)=1\}|\le 1$
for all $x\in\Z^d$. This implies that at most $\sw(\A^{x,t})$
cylinders of $C(\zeta)$ can be incompatible with cylinders in
$\A^{x,t}$. Hence, \reff{2289} is bounded by
\begin{equation}
  \label{231}
 1-(1-e^{-(1-b)t})^{\sw(\A^{x,t})}\,.
\end{equation}
The expectation of \reff{231} is given by
\begin{eqnarray}
  \label{233}
 \sum_{n\ge 1} \Bigl[1-(1-e^{-(1-b)t})^n\Bigr]\,
\P\Bigl(\sw(\A^{x,0})=n\Bigr)
\end{eqnarray}
because the distribution of $\A^{x,t}$ does not depend on $t$. Our
hypothesis of no backwards oriented percolation implies that
$\A^{x,0}$ contains a finite number of (finite) contours. Hence
$\sum_{n\ge 1} \P(\sw(\A^{x,0})=n)=1$ and by dominated convergence
\reff{233} goes to zero as $t\to\infty$. This proves \reff{754}.

To prove \reff{335} we start from the expectation of \reff{228} and use
\reff{233} to bound the expected value of the second term by
\begin{eqnarray}
   \lefteqn{e^{-(1-b)t} \sum_{n\ge 1} \P\Bigl(\sw(\A^{x,0})=n\Bigr)\, 
\sum_{k=0}^{n-1} \Bigl(1-e^{-(1-b)t}\Bigr)^{k} } \nonumber\\
  & \le & e^{-(1-b)t} \sum_{n\ge 1} n\, \P\Bigl(\sw(\A^{x,0})=n\Bigr)
  \nonumber \\
  & \le & e^{-(1-b)t}\,\, \E\Bigl(\sw(\A^{x,0})\Bigr). \label{31.a}
\end{eqnarray}

The above arguments prove that the process converges, uniformly in the initial
configuration, to the invariant measure $\mu$. An immediate consequence is
that $\mu$ is the unique invariant measure. This concludes the proof of (1)
and (2) of Theorem \ref{225}.\ \square

\subsection{Finite-volume effects}\label{s131}

To prove inequality \reff{73a} we use \reff{505}, \reff{506} and
\reff{77} to get 
\begin{eqnarray}
  \mu f- \mu^\Lambda f\;=\; \E f(\eta_0)- \E f(\eta^\Lambda_0)\; =\; \E
  \Bigl[\Phi\Bigl(\A(\supp(f))\Bigr)- 
  \Phi\Bigl( \A^{\Lambda}(\supp(f))\Bigr)\Bigr],\label{81a}
\end{eqnarray}
where $\Phi$ is the function referred to in \reff{77a}.
By definition, 
\begin{equation}
 \Phi\Bigl(\A(\supp(f))\Bigr) \;\le\; \|f\|_\infty\;. 
\label{107}
\end{equation}
Hence, inequality \reff{73a} follows from \reff{81a}.

Since the spatial projections of the set of ancestors of $\supp(f)$ are
finite, the right hand side of \reff{73a} goes to zero, proving, in
particular, the weak convergence of $\mu^\Lambda$ to $\mu$. \square

\subsection{Mixing.  Its relation with a coupling construction}

The proof of (4) of Theorem \ref{225} is very similar in spirit to the
above proof but it requires a somewhat more delicate argument based on
the coupling of two continuous-time versions of the backwards
percolation process.  We first notice that \reff{73d} is a
straightforward consequence of \reff{73b}, because in the absence of
backwards percolation, the spatial projections of the set of ancestors
of $\supp(f)$ and $\supp(g)$ are finite. This implies that the right
hand side of \reff{73b} goes to zero.

To prove \reff{73b} we use \reff{505} and \reff{77} to get
\begin{eqnarray}
  \lefteqn{\vert \mu (f g) - \mu f\, \mu g \vert}\nonumber\\
  &=&  \E(f(\eta_0)g(\eta_0))- \E f(\eta_0)\,\E g(\eta_0)\nonumber\\
&=& \E \Bigl[\Phi\Bigl(\A(\supp(f))\Bigr)\,
 \Phi\Bigl(\A(\supp(g))\Bigr)\Bigr]-\E\Bigl[
  \Phi\Bigl(\widehat \A(\supp(f))\Bigr)\,\Phi\Bigl(\widehat
  \A(\supp(g))\Bigr)\Bigr]\;, \nonumber\\
\ \label{81}
\end{eqnarray}
where $\Phi$ is the function referred to in \reff{77} and
$(\widehat\A(\supp(f)),\widehat\A(\supp(g)))$ has the \emph{same}
marginal distributions as $(\A(\supp(f)),\A(\supp(g)))$ but its
marginals are \emph{independent}.

Identity \reff{81} shows that to obtain \reff{73d} it is enought to
construct a coupling (joint construction) of the four processes
$$
\Bigl(\A(\supp(f))\,,\,\A(\supp(g))\,,\,\widehat\A(\supp(f))\,,\,
\widehat\A(\supp(g))\Bigr)\;,
$$ 
such that
\begin{equation}
\label{555}
\A(\supp(f))=\widehat\A(\supp(f))
\end{equation}
and
\begin{equation}
\label{556}
\widehat\A(\supp(g))\sim\A(\supp(f))\quad\hbox{ implies }\quad
\A(\supp(g))=\widehat\A(\supp(g) )\;.
\end{equation}
Indeed, from \reff{107} and \reff{555}--\reff{556} we obtain that the
last line of \reff{81} is bounded above by the right-hand side of
\reff{73b}. 

In the remaining of the section we discuss the construction of the 
coupling with properties \reff{555}--\reff{556}.  The construction is
natural and straightforward, but unavoidably technical.  As an
alternative we mention the approach based on ``duplicated variables''
(von Dreifus, Klein and Perez, 1995; Bricmont and Kupiainen, 1996),
which is probabilistically simpler but requires some combinatorial
input.  

\subsection{Construction of a four-clan coupling} \label{ssfc}

We need to couple two clans in the same random set of cylinders with
two independent copies with the same marginal distribution.  Moreover,
to strengthen our results we need to ensure that the marginal
realizations remain the same as much as possible.  The coupling
(Section \ref{S30}), is based on a construction of backwards
percolation clans as non-homogeneous continuous time Markov
processes (Section \ref{118}).  The hypothesis on the absence of
non-oriented percolation for some time interval $(0,h)$ is needed
for the infinite-volume construction of the coupling.

\subsubsection {A continuous-time construction of the backwards percolation
  clan} \label{118}

For $\Upsilon\subset
\Z^d$ define
\begin{equation}
  \label{65}
  \A_t(\Upsilon)\; =\; \Bigl\{C'\in \A(\Upsilon):
  0>\birth(C')>-t\Bigr\}\; =\; \A(\Upsilon)\cap \C[-t,0] ,
\end{equation}
that is, the set of cylinders in $A(\Upsilon)$ with birth-time
posterior to $-t$.  The inclusion of a new cylinder in the time
interval $[t,t+h]$ depends on the existence of a birth Poisson mark in
$[-t-h,-t]$ whose corresponding cylinder is incompatible with some
$C'\in \A_t(\Upsilon)$.  That is, if $C$ is a cylinder with 
$\basis(C)\not\sim \basis(C')$ for some $C'\in \A_t$, 
\begin{eqnarray}
\lefteqn{\P \Bigl(\A_{t+h}=\widetilde\A\cup
 C\Given\A_t=\widetilde \A,\, \A_{t'}=\widetilde
 \A_{t'}, \,t'\in[0,t) \Bigr) }\nonumber\\ 
 &=& \P \Bigl\{C\in\C\,:\,\birth(C)\in [-t-h,-t]\,,\,
 \death(C)>t-\ti(\widetilde\A,\basis(C))\Bigr\}\;+ \;o(h)   \;.
 \nonumber \\
\ 
\end{eqnarray}
We have denoted 
\begin{equation}
\ti(\A_t,\ga) = \min\{\birth(C'):C'\in \A_t,
\basis(C')\not\sim\ga\}
\end{equation}
and abbreviated $\A_t(\Upsilon)=\A_t$.  The remainder $o(h)$ is the
correction related to the probability that $C$ is not the only
cylinder born in $[-t-h,-t]$. Since the birth-time is independent of
the lifetime which is exponentially distributed with rate one,
\begin{eqnarray}
 \lefteqn{\P \Bigl(\A_{t+h}=\widetilde\A\cup
 C\Given\A_t=\widetilde \A ,\, \A_{t'}=\widetilde
 \A_{t'}, \,t'\in[0,t)  \Bigr) }\nonumber\\ 
 &=& \P \Bigl\{C\in\C\,: \,\birth(C)\in [-t-h,-t]\Bigr\}\,
\P\Bigl( \life(C)>t-\ti(\widetilde\A,\basis(C))\Bigr)\;
 +\; o(h)  \nonumber\\
 &=& h\,e^{-\beta|\basis(C)|}\, e^{-t+\ti(\widetilde\A,\basis(C))} 
\;+\; o(h)\;.\label{oo2}
\end{eqnarray}
This implies that when the configuration at time $t^-$ is
$\widetilde\A$, a new cylinder with basis $\gamma$ is included in
$\A_t(\Upsilon)$ at rate
\begin{equation}
  \label{ooo}
  e^{-\beta|\gamma|}\,\, e^{-t+\ti(\widetilde\A,\gamma)}.
\end{equation}
From \reff{oo2}, as in the computation of the forward Kolmogorov
equations, we get
\begin{equation}
  \label{oo1}
  \E\,\Bigl({d f(\A_t)\over dt}\Given \A_s,\, 0\le s\le t\Bigr) 
  \;=\; \sum_\gamma \int_{t-\ti(\A_t,\ga)}^\infty 
  ds\,\,e^{-s}\, e^{-\beta|\gamma|}\,
   \Bigl[f(\A_t\cup (\ga,t,s)) - f(\A_t)\Bigr]
\end{equation}
where the sum is over the set $\{\gamma\in\G: \gamma\not\sim
\basis(C')$ for some $C'\in\A_t\}$.
This equation characterizes the law of the process $\A_t(\Upsilon)$ as
a non-homogeneous Markov process.

We now construct $\A_t(\Upsilon)$ by combing the Poisson marks backwards
in time in a continuous manner. 

\paragraph{Finite-volume case.}
If we only consider contours contained in a finite set
$\Lambda\subset\Z^d$, there is only a finite set of possible bases for
the cylinders and the Poisson marks are well ordered with probability
one.  The construction proceeds mark by mark backwards in time. Set
$\A_0(\Upsilon) = \Upsilon $. If there is a Poisson birth mark at time
$-t$ whose corresponding cylinder is called $C''$, then
\begin{itemize}
\item if $C''\not\sim C'$ for some $C'\in \A_{t-}(\Upsilon)$, set
  $\A_t(\Upsilon)= \A_{t-}(\Upsilon)\cup\{C''\}$.
\item if $C''\sim C'$ for all $C'\in \A_{t-}(\Upsilon)$, set
  $\A_t(\Upsilon)= \A_{t-}(\Upsilon)$.
\end{itemize}
where the incompatibility between cylinders was defined in \reff{502}.

\paragraph{Infinite volume case.}
In infinite volume the construction can be performed using a
percolation argument as in Section \ref{S34}.  By hypothesis, there
exists an $h$ such that each cylinder born in the interval $[-h,0]$
belongs to a finite non-oriented clan. Hence the set of cylinders born
in the interval $[-h,0]$ can be partitioned in connected families:
\begin{equation}
 \label{hk}
 \C[-h,0] = \bigcup_{k\ge 0}^{\cdot} \H_k[-h,0]
\end{equation}
where the sets $\H_k[s,t]$ are the maximal sets of cylinders with the
property that cylinders in different $\H_k$'s are compatible.  We can
then well-order the birth-time of the cylinders inside each $\H_k$ and
proceed as for the finite-volume case.  This yields the process
$\A_t(\Upsilon)$ for $t\in[0,h]$. To extend the construction for
arbitrary $t>0$, we simply repeat the previous procedure in $[h,2h]$,
$[2h,3h]$, etc.

\subsubsection {A coupling between two interacting and two independent clans}
\label{S30}

We take two independent marked Poisson process whose marks and
cylinders we respectively call blue and red. We enlarge
our probability space and continue using $\P$ and $\E$ for the
probability and expectation with respect to the space
generated by the product of the blue and red Poisson processes.  Using
these marks we construct simultaneously the processes
$(\A_t(\Lambda_1),\A_t(\Lambda_2),\widehat\A_t(\Lambda_1),
\widehat\A_t(\Lambda_2))$, for $\Lambda_1,\Lambda_2 \subset \Z^d$,
in the following way. 
\begin{enumerate}
\item The processes $\A_t(\Lambda_1)$ and $\A_t(\Lambda_2)$ are constructed
using only the blue marks, as described in Subsection \ref{118}, and ignoring
the red marks. Hence, they are the clans of $\Lambda_1$ and
$\Lambda_2$ respectively. 
\item The process $\widehat\A_t(\Lambda_1)$ is also constructed only
  with the blue marks, hence it coincides with $\A_t(\Lambda_1)$. 
\item \label{pe} The process $\widehat\A_t(\Lambda_2)$ is constructed with a
  precise combination of blue and red marks in such a way that (a) it
  coincides with $\A_t(\Lambda_2)$ for a time interval that is as
  long as possible; (b) it is independent of
  $\widehat\A_t(\Lambda_1)$, and (c) it has the same marginal
  distribution as $\A_t(\Lambda_2)$.
\end{enumerate}
Property \reff{pe} is achieved in the following way.

\paragraph{Finite-volume case.}
If both $\Lambda_1$ and $\Lambda_2$ are finite sets, we order the
marks by appearance and introduce a flag variable, Flag$(t) \in
\{0,1\}$, which indicates if some cylinder of $\widehat\A_t(\Lambda_1)$
is incompatible with some cylinder of $\widehat\A_t(\Lambda_2)$:
\begin{equation}
  \label{71}
  \hbox{Flag}(t)
\;=\; \one\Bigl\{C'\not\sim C''\hbox{ for some
  }C'\in\widehat\A_t(\Lambda_1),C''\in\widehat\A_t(\Lambda_2)\Bigr\}.
\end{equation}
We now proceed as follows, mark by mark backwards in time.  First, we
set Flag$(0)=0$. The construction guarantees that Flag$(t)=0$ implies
$\A_t(\Lambda_i) = \widehat\A_t(\Lambda_i)$ for $i=1,2$.
\begin{itemize}
\item If at time $-t$ a (blue or red) mark is present and  Flag$(t-)=0$, then
\begin{itemize}
\item If the mark is blue and the corresponding cylinder can be included in
  $\A_{t-}(\Lambda_1)$ but not in $\A_t(\Lambda_2)$, then include it in
  $\A_t(\Lambda_1)$ and $\widehat\A_t(\Lambda_1)$. Analogously, if it can be
  included in $\A_{t-}(\Lambda_2)$ but not in $\A_t(\Lambda_1)$, then include
  it in $\A_t(\Lambda_2)$ and $\widehat\A_t(\Lambda_2)$. Keep the Flag
  $=0$.
 
\item If the mark is blue and the corresponding cylinder can be included in
  both $\A_{t-}(\Lambda_1)$ and $\A_{t-}(\Lambda_2)$, then include it in both
  $\A_{t}(\Lambda_1)$ and $\A_{t}(\Lambda_2)$ but include it only in
  $\widehat\A_t(\Lambda_1)$. Set the Flag $=1$.
 
\item If the mark is red and the corresponding cylinder can be included in
  both $\A_{t-}(\Lambda_1)$ and $\A_{t-}(\Lambda_2)$, then include it only in
  $\widehat\A_t(\Lambda_2)$. Set the Flag $=1$.
 
\item If the mark is red and the corresponding cylinder can be
  included in either $\A_{t-}(\Lambda_1)$ or $\A_{t-}(\Lambda_2)$
  but not in both of them, then ignore the mark. Keep the Flag $=0$.
\item If the mark is red or blue but the corresponding cylinder can be
  included in neither $\A_{t-}(\Lambda_1)$ nor $\A_{t-}(\Lambda_2)$ then
  ignore the mark. Keep the Flag
  $=0$.
\end{itemize}
\item If at time $-t$ a (blue or red) mark appears and Flag$(t-)=1$, then use
  blue marks for $\A_t(\Lambda_1)$, $\A_t(\Lambda_2)$, and
  $\widehat\A_t(\Lambda_1)$ and red marks for $\widehat\A_t(\Lambda_2)$.
  Keep the Flag $=1$.
\end{itemize}

To verify that the above coupling has the right marginals, it
suffices to notice that the rate of inclusion of cylinders in each one
of the marginals is precisely given by \reff{ooo}. 

\paragraph{Infinite-volume case.}
In the case in which at least one of $\Lambda_1$ and $\Lambda_2$ is
infinite, we consider families $\overline{\H}_k$ analogous to the
$\H_k$ given in \reff{hk} but defined using the time interval
$[-h/2,0]$ and both red and blue marks. Therefore in the combined set
of cylinders there is no non-oriented percolation and we can construct
the coupling working in a finite set of $\overline{\H}_k$'s at a time.
We then continue working in time intervals of lenght $h/2$ to reach
arbitrary times $t$.

By construction, the coupling satisfies \reff{555}. It also satisfies
\reff{556} because the flag changes from $0$ to $1$ (and remains $1$
forever) the first time a cylinder of $ \A_t(\supp(g))$ is
incompatible with a cylinder of $\widehat \A_t(\supp(f))$. This
implies
\begin{equation}
  \label{85}
  \hbox{Flag}(\infty)=\one\Bigl\{\A(\supp(f))\not\sim \widehat
  \A(\supp(g))\Bigr\}.
\end{equation}

\section{Branching processes. Time length and space width}\label{SSS}

In this section we estimate the time-length and space-width of the
families of ancestors $\A^{x,t}$.  We follow the well known approach
of introducing a branching process that dominates the backward
percolation process (see eg.\ Hall, 1985), though we must consider
\emph{multitype} branching.  The main result of this section is
Theorem \ref{thm:hall} which shows that the hypotheses of Theorem
\ref{18a} lead to exponential upper bounds of both $\tl(\A^{x,t})$
and $\sw(\A^{x,t})$.

\subsection{Multitype branching processes} \label{mbp}

We introduce a multitype branching process $\B_n$, in the set of
cylinders, which dominates $\A_n$.  To do this we look ``backwards in
time'' and let ``ancestors'' play the role of ``branches''.  In
particular, births in the original marked Poisson process correspond
to dissapearance of branches.  We reserve the words ``birth'' and
``death'' for the original forward-time Poisson process.  

We start by enlarging our probability space and defining, for any
given set $\{C_1,\dots,C_k\}$, \emph{independent} random sets
$\B_1^{C_i}$ with the same marginal distribution as $\A_1^{C_i}$.  The
important point here is that
\begin{equation}
  \label{wt}
  \bigcup_{i=1}^{k}\A_1^{C_i} \subset \bigcup_{i=1}^{k}\B_1^{C_i}.
\end{equation}
The proof of this fact relies in fixing a way to distribute common
ancestors.  For example, consider the total order $\prec$ in the set of
cylinders induced by the birth-times. That is $C \prec C'$ if and only
if $\birth(C)\leq \birth(C')$.  For any finite set of cylinders
$\{C_1,\dots,C_k\}$ such that $C_i\prec C_{i+1}$, $i=1,\dots,k-1$,
define
\begin{equation}
  \label{ta}
  \widetilde\A_1^{C_j}
= \A_1^{C_j} \setminus \left(\bigcup_{\ell = 1}^{j-1} \A_1^{C_\ell}\right).
\end{equation}
This ensures that $\widetilde\A_1^{C_j}$ are independent sets and
\begin{equation}
  \label{cy}
  \bigcup_{i=1}^{k}\A_1^{C_i} = \bigcup_{i=1}^{k}\widetilde\A_1^{C_i}.
\end{equation}
On the other hand,
for any $C$, $\widetilde\A_1^C$ is stochastically dominated by $\A_1^C$ [that
is, there exists a joint realization $(\widetilde\A_1^C, \A_1^C)$ such that
$\P(\widetilde\A_1^C\subset\A_1^C) = 1$]. From this observation and
\reff{cy} we get \reff{wt}.

The procedure defined by $\B_1$ naturally induces a multitype branching process
in the space of cylinders.  We define the $n$-th generation of the branching
process by
\begin{equation}
  \label{an1}
  \B_n^C = \{ \B_1^{C'}: C'\in \B_{n-1}^C\}
\end{equation}
where for all $C'$, $\B_1^{C'}$ has the same distribution as $\A_1^{C'}$ and are
independent random sets depending only on $C'$. Inductively,
\begin{equation}
  \label{4b}
  \A_{n}^C\subset \B_{n}^C.
\end{equation}
Indeed, 
\begin{equation}
  \label{4a}
  \A_n^C = \bigcup_{C'\in \A_{n-1}^C} \A_1^{C'} =  \bigcup_{C'\in
  \A_{n-1}^C} \widetilde\A_1^{C'} 
\end{equation}
where in the definition of $\widetilde\A_1^{C'}$ we use $\{C_1, \dots, C_k\}
=\cup_{i=0}^{n-1}\A_{i}^C$. Hence, the inductive hypothesis
$\A_{i}^C\subset \B_{i}^C$, for $i=1,\dots,n-1$, yields \reff{4b}.

In consistence with our previous notation, we denote
\begin{equation}
  \label{608}
  \B^C =\bigcup_{n\ge 0} \B_n^C\;\;;\;\;\;\; 
\B^{x,t} =\bigcup_{n\ge 0} \B_n^{x,t}\;\;;\;\;\;\; 
\B^{\Upsilon} =\bigcup_{x\in\Upsilon} \B^{x,0},
\end{equation}
the branching clans of $C$, $(x,t)$ and $\Upsilon$ (at time $0$)
respectively. 
By \reff{4b},
\begin{equation}
  \label{610}
  \A^C \subset \B^C ;\;\; \A^{x,t}\subset\B^{x,t} ;\;\;
  \A^{\Upsilon}\subset\B^{\Upsilon}.
\end{equation}
Defining the time length and space width of this clan as in \reff{229} and
\reff{232}, we get
\begin{equation}
  \label{611}
  \tl(\A^C)\le \tl(\B^C) \;;\; \tl(\A^{x,t})\le\tl(\B^{x,t}) \;;\;
  \tl(\A^{\Upsilon})\le\tl(\B^{\Upsilon}) \;,
\end{equation}
and similarly for the respective space widths.

The (multitype) branching process $\B_n$ induces naturally a multitype
branching process in the set of contours. For a cylinder $C$ with
basis $\ga$ and birth-time $0$, define $\b^\ga_n\in \N^\G$ as the
number of cylinders in the $n$th generation of ancestors of $C$ with
basis $\th$\/:
\begin{equation}
  \label{10}
 \b^\ga_n(\th) = \Bigl|\Bigl\{C'\in \B_n^C:\basis(C') = \th \Bigr\}
 \Bigr|\;.
\end{equation}
This process will be useful in estimating the space properties of the
clans of ancestors.  We have the following relationship:
\begin{equation}
  \label{14}
  \sum_\th \b^\ga_n(\th) = |\B_n^C|.
\end{equation}

The process $\b_n$ is a multitype branching process whose offspring
distributions are Poisson with means
\begin{eqnarray}
  \mm(\ga,\th) & = & \one\{\ga \not\sim \th \}\ e^{-\beta |\th|}
  \int_{0}^{\infty} e^{-t} dt \nonumber \\
  & = & \one\{\ga \not\sim \th \} e^{-\beta |\th|}.
\label{eq:mu}
\end{eqnarray}
To see this, notice that the cylinders $C'$ with basis $\theta$ that
are potential ancestors of $C$ (with basis $\ga$) form a Poisson
process of rate $\one\{\ga \not\sim \th \}\,e^{-\beta |\th|}$. Each of
those cylinders is an ancestor of $C$ if its lifetime is bigger than
the difference between the birth-time of $C$ and $C'$.  The lifetimes
of different cylinders are independent exponentially distributed
random variables of rate $1$. The probability that the lifetime of any
given cylinder is bigger than $t$ is given by $e^{-t}$. Hence, the
birth-times of the ancestors of $C$ with basis $\th$ form a (non
homogeneous) Poisson process of rate depending on $t$ given by
$\one\{\ga \not\sim \th \}\,e^{-\beta |\th|}\,e^{-t}$. The mean number
of births is therefore given by \reff{eq:mu}.

\begin{lemma}
  \label{16k} The means \reff{eq:mu} satisfy
\begin{equation}
  \label{16}
  \sum_{\th} \mm^{n}(\ga,\th) \;\le\;
  \sum_{\th} |\th|\,\mm^{n}(\ga,\th) \;\leq\; \vert \ga \vert\, \alpha^n,
\end{equation}
where $\alpha$ is defined in \reff{19}.
\end{lemma}
\proof
\begin{eqnarray}
  \label{a16}
  \sum_{\th} |\th|\, \mm^{n}(\ga,\th) &=&
    \sum_{\ga_1:\ga_1\not\sim\ga} e^{-\be\vert 
    \ga_1\vert}\sum_{\ga_2:\ga_2\not\sim\ga_1} e^{-\be\vert \ga_2\vert}\ldots
  \sum_{\th:\th\not\sim\ga_{n-1}} |\th|\,
  e^{-\be\vert \th\vert} \nonumber\\[10pt]
  &=& |\ga| \sum_{\ga_1:\ga_1\not\sim\ga} {|\ga_1|\over |\ga|} e^{-\be\vert
    \ga_1\vert}\sum_{\ga_2:\ga_2\not\sim\ga_1} {|\ga_2|\over |\ga_1|}
  e^{-\be\vert \ga_2\vert}\ldots \sum_{\th:\th\not\sim\ga_{n-1}}
  {|\th|\over |\ga_{n-1}|}   e^{-\be\vert \th\vert} \nonumber\\[10pt]
  &\leq& \vert \ga \vert \biggl(\sup_\ga \sum_{\th:\th\not\sim\ga} {|\th|\over
    |\ga|}e^{-\be\vert \th\vert}\biggr)^n .\ \square
\end{eqnarray}

This lemma shows, in particular, that the branching process $\b_n$ is
subcritical if $\alpha<1$.  

\subsection{Continuous-time branching process} \label{ctbp}

Let $C$ be a cylinder with basis $\ga$ and birth-time $0$. Combing backwards
continuously in time the branching clan $\B^C$ we define a continuous-time
multitype branching process $\psi^\ga_t(\th)= $ number of contours of type
$\th$ present at time $t$ (of this process) whose initial configuration is
$\de_\ga$. Each $C'\in\B^C$ is a branch, that is, belongs to the first
generation of ancestors of a unique cylinder $U(C')$ in $\B^C$.  In
the branching process $\psi_t$ all the branches (ancestors) of $U(C')$
appear simultaneously at the birth of $U(C')$, that is when $U(C')$
dissapears if we look backwards in time.  Therefore the part of $C'$
in the interval $[\birth(U(C')),\death(C')]$ is ignored.  Formally, 
\begin{equation}
  \label{25}
  \psi^\ga_t(\th)= \Bigl|\Bigl\{C'\in \bar\B^C: \basis(C')=\th, \birth(C')<-t<
  \birth(U(C')), \life(C')\ni t\Bigr\}\Bigr|.
\end{equation}

In the process $\psi_t$, each contour $\ga$ lives a mean-one exponential
time after which it dies and gives birth to $k_\th$ contours $\th$,
$\th\in \G$, with probability
\begin{equation}
  \label{30}
  \prod_\th {e^{\mm(\ga,\th)} \,\mm(\ga,\th)^{k_\th} \over k_\th!}
\end{equation}
for $k_\th\ge 0$. These are independent Poisson distributions of mean
$\mm(\ga,\th)$. The infinitesimal generator of the process is given by
\begin{equation}
  \label{30.1}
  L f (\psi) = \sum_{\ga \in \G} \psi(\ga) \sum_{\eta \in {\cal
  Y}_{0}(\ga)} \prod_{\th: \eta(\th) \ge 1} {e^{\mm(\ga,\th)}\,
  \mm(\ga,\th)^{\eta(\th)} \over \eta(\th)!}\, \left[ f(\psi
  + \eta - \delta_{\ga}) - f(\psi) \right]
\end{equation}
where $\psi, \eta \in {\cal Y}_{0} = \{ \psi \in \N^\G;
\sum_\th \psi(\th) < \infty\}$ and ${\cal Y}_{0}(\ga) = \{ \psi \in {\cal
Y}_0; \psi(\th) \ge 1 \mbox{ implies } \th \not\sim \ga \}$ and
$f: {\cal Y}_0 \rightarrow \N$.

The branching process $\psi^\ga_t$ allows us to estimate the
time-length of a clan, due to the obvious fact:
\begin{equation}
  \label{25a}
  \sum_\th\psi^\ga_t(\th) = 0\, \hbox{ implies }\, \tl(\B^C)<t\;. 
\end{equation}
Let $M_t(\ga,\th)$ be the mean number of contours of type $\th$ in $\psi_t$ and
$R_t(\ga)$ its sum over $\th$:
\begin{equation}
  \label{28}
  M_t(\ga,\th)= \E \psi^\ga_t(\th)\ ;\ \ \ \ R_t(\ga)=\sum_\th M_t(\ga,\th)\;.
\end{equation}
The bound we need is given in the next lemma.

\begin{lemma} The mean number of branches, $R_t(\ga)$ satisfies
  \label{33}
  \begin{equation}
    \label{34}
   \P\Bigl(\sum_\th \psi_t^\ga(\th) >0\Bigr)\, \le \,R_t(\ga)\,\le\, |\ga|\,
   e^{(\alpha-1)t}. 
  \end{equation}
\end{lemma}
\proof The first inequality is immediate because $\sum_\theta\psi_t^\ga(\th)$
assumes non-negative integer values and $R_t(\ga)$ is its mean value.

To show the second inequality we first use the generator given by
\reff{30.1} to get the Kolmogorov backwards equations for $R_t(\ga)$:
\begin{equation}
  \label{27}
  {d\over dt} R_t(\ga) = \sum_{\ga'} \mm(\ga,\ga')
  R_t(\ga') - R_t(\ga).
\end{equation}
Since $R_0(\ga')\equiv 1$, the solution is
\begin{equation}
  \label{29}
  R_t(\ga) = \sum_{\ga'}\biggl[ \exp [t( \mm-I)]\biggr](\ga,\ga')
\end{equation}
where $\mm$ is the matrix with entries $\mm(\ga,\ga')$ and $I$ is the
identity matrix.  This can be rewritten as
\begin{eqnarray}
  \label{344}
  R_t(\ga)&= & e^{-t}\sum_{n\ge 0} {t^n \over n!}
  \sum_{\ga'} \mm^n(\ga,\ga') \\[5pt]
  &\le& e^{-t}\sum_{n\ge 0} {t^n \over n!}  \, |\ga| \, \alpha^n \;,
\end{eqnarray}
where the last bound is just the leftmost inequality in \reff{16}. \square

\subsection{Time length and space width}\label{5.4}

We are now ready to provide bounds for the time length and space width
of the percolation clan.

\begin{thm} \label{thm:hall} If $\beta>\beta^*$ (i.e.\/ $\alpha(\beta)< 1$), then
  \begin{itemize}
  \item[(i)] The probability of backward oriented percolation is zero.
\item[(ii)] For any positive $b$,
  \begin{equation}
    \label{336}
    \P\Bigl( \tl(\A^{x,t})>bt\Bigr) \;\le\; \alpha_0\, e^{-(1-\alpha)bt}
  \end{equation}

\item[(iii)]
  \begin{equation}
    \label{46a}
\E\Bigl( \sw(\A^{x,t})\Bigr) \;\le\; \frac{\alpha_0(\beta)}{1-\alpha(\beta)}
  \end{equation}

\item[(iv)]
  \begin{equation}
    \label{46b}
\E\Bigl( \exp[a\,\sw(\A^{x,t})]\Bigr)\; \le\;
\frac{\alpha_0(\beta-a)}{1-\alpha(\beta-a)} 
  \end{equation}

\item[(v)]
  \begin{equation}
    \label{46c}
    \P\Bigl( \sw(\A^{x,t})\ge \ell\Bigr) \;\le\;
    \frac{\alpha_0(\widetilde\beta)}{1-\alpha(\widetilde\beta)}\,
    e^{-(\beta-\widetilde\beta)\,\ell}    
  \end{equation} 
for any $\widetilde\be\in(\beta^*,\beta)$.
 \end{itemize}
\end{thm}

\paragraph{Proof} (i)  We follow an idea of Hall
  (1985). For each $C \in \C$ we use the domination \reff{4b} and the
  identity \reff{14}.  Therefore, to prove that there is no backward oriented
  percolation it is enough to prove that, for fixed $\ga$
\begin{equation}
  \P \biggl( \sum_{\th} \b^\ga_{n}(\th) \neq 0 \mbox{ for infinitely
    many $n$} \biggr) = 0.
\label{eq:6}
\end{equation}
Since $\b^\ga_n(\th)$ assumes non negative integer values, 
by Borel-Cantelli lemma, a sufficient condition for
(\ref{eq:6}) is
\begin{equation}
\sum_{n} \sum_{\th} \mm^{n}(\ga,\th)<\infty. \label{eq:5}
\end{equation}
But this follows from Lemma \ref{16k}. 
\smallskip

(ii) By \reff{25a} and \reff{611} for each $x \in \ga$
and $s \le t$
\begin{equation}
  \label{31}
  \sum_\th \psi^\ga_{t}(\th) = 0 \hbox{\ \ implies\ \ } \tl(\A^{x,0}) \le t.
\end{equation}
Hence,
\begin{equation}
  \label{339}
  \P(\tl(\A^{x,0})>t) \;\le\; \sum_{\ga\ni x} \P(\eta_0(\ga)=1)\, R_t(\ga)
  \le\; \sum_{\ga\ni x} e^{-\beta|\ga|}\, R_t(\ga)
  \;\le\;  \alpha_0\, e^{-(1-\alpha)t}
\end{equation}
\smallskip
by the rightmost inequality in \reff{34}.

(iii) We find
upperbounds for the space diameter of the backwards percolation clan
through upperbounds for the total number of  
occupied points by the multitype branching process $\b_n$ defined by
\reff{10}. In fact,
\begin{equation}
 \label{33.a}
 \sw(\A^{x,0}) \le \sum_{\ga\ni x}\eta_0(\ga)\sum_{n} \sum_{\th} |\th|\,
 \b^{\ga}_n(\th).
\end{equation}
By \reff{623}, $\E\eta_0(\ga)\le e^{-\beta|\ga|}$, hence by \reff{16}
\begin{eqnarray}
  \E \Bigl(\sum_{\ga\ni x}\eta_0(\ga) \sum_{n} \sum_{\th} |\th|\,
  \b^{\ga}_n(\th) \Bigr) 
& \le & \sum_{\ga\ni x}e^{-\beta|\ga|} \sum_{n} \sum_{\th} |\th|\,
  \mm^{n}(\ga,\th) \nonumber \\
  & \le & \alpha_0 \sum_{n} \alpha^n. 
\end{eqnarray}

(iv) Write
\begin{eqnarray}
   \E(e^{a\sw}) &=& \sum_{\ell} e^{a\ell}\, \P(\sw=\ell) \nonumber\\
&\le& \sum_{\ell}e^{a\ell}\sum_k \sum_{\ga_1,\dots,\ga_k}
\one\Bigl\{|\ga_1\cup\dots\cup\ga_k| = \ell\Bigr\}\, 
\P\Bigl(\ga_1\ni 0, \b_1^{\ga_1}(\ga_2)\ge 1,\dots,
\b_1^{\ga_{k-1}}(\ga_k)\ge 1\Bigr)\;.\nonumber\\ 
 \label{700}
\end{eqnarray}
By the Markovian property of $\b_n$ we get
\begin{eqnarray}
  \lefteqn{ \P\Bigl(\ga_1\ni 0,
    \b_1^{\ga_1}(\ga_2)\ge 1,\dots, \b_1^{\ga_{k-1}}(\ga_k)\ge
    1\Bigr)}\nonumber\\
  &=& \P(\ga_1\ni 0)\,\P(\b_1^{\ga_1}(\ga_2)\ge 1)\cdots
  \P(\b_1^{\ga_{k-1}}(\ga_k)\ge
  1)\nonumber\\
  &\le& \one\Bigl\{\ga_1\ni 0, \ga_1\not\sim\ga_2, \dots,
    \ga_{k-1}\not\sim\ga_k \Bigr\}\,\prod_{i=1}^k e^{-\be|\ga_i|}
    \label{701} 
\end{eqnarray}
Substituting this in \reff{700} and using that 
\begin{equation}
  \label{703}
 e^{a\ell}\, \one\Bigl\{|\ga_1\cup\dots\cup\ga_k| = \ell\Bigr\}
\;\le\; \one\Bigl\{|\ga_1\cup\dots\cup\ga_k| = \ell\Bigr\} \, 
\exp\Bigl(a\sum_{i=1}^k |\ga_i|\Bigr)
\end{equation}
we get
\begin{eqnarray}
  \label{702}
  \E(e^{a\sw})&\le& \sum_k \sum_{\ga_1,\dots,\ga_k} \one\Bigl\{\ga_1\ni 0,
  \ga_1\not\sim\ga_2, \dots, \ga_{k-1}\not\sim\ga_k \Bigr\}\,
  \exp\Bigl(-(\be-a)\sum_{i=1}^k |\ga_i|\Bigr)\nonumber\\
  &=& \sum_k \sum_{\ga_1\ni\,
    0}|\ga_1|\,e^{-(\be-a)|\ga_1|}{1\over|\ga_1|}\sum_{\ga_2\not\sim\ga_1}|\ga_2|
  e^{-(\be-a)|\ga_2|}\dots {1\over|\ga_{k-1}|} \sum_{\ga_k\not\sim\ga_{k-1}}
  e^{-(\be-a)|\ga_{k}|}\nonumber\\
  &\le& \alpha_0(\be-a) \sum_{k\ge 0} \alpha(\be-a)^k.\ 
\end{eqnarray}

(v) It suffices to use (iv) and the exponential Chevichev inequality and to
notice that $a$ must be less than $\be - \be^*$ to avoid a zero in the
denominator of \reff{46b}.  \square

\medskip
\noindent{\bf Remarks. } {\bf 1.} 
Part (ii) can in fact be proven by a more elementary argument not
requiring the continuous-time construction of Section \ref{ctbp}. The
argument gives the same rate of decay as in \reff{336} but a worse
leading constant. Let us sketch it.
\begin{equation}
  \label{800}
  \P\Bigl(\tl(\A^{x,0})>t\Bigr) \;\le\; \sum_k
  \sum_{\ga_1,\dots,\ga_k} \one\{\ga_1\ni 0, \ga_1\not\sim\ga_2,
  \dots, \ga_{k-1}\not\sim\ga_k \}\, \P(S_1+\cdots+S_k>t)
\end{equation}
where $S_i$ are independent mean one exponentially distributed random
variables and independent of $\ga_i$. The time $S_i$ represents the period
between the birth of $\ga_i$ and $\ga_{i+1}$.  As the sum of independent
exponentials is a gamma distribution with parameters $k$ and $1$, 
\begin{equation}
  \label{801}
  \P(S_1+\cdots+S_k>t) \;=\; e^{-t}\, \sum_{i=0}^k {t^i\over i!}\;.
\end{equation}
Therefore \reff{800} is bounded by
\begin{eqnarray}
 &&\hspace{-2cm} e^{-t}\,\sum_{i=0}^\infty {t^i \over i!} \sum_{k\ge i}
    \sum_{\ga_1\ni\,
      0}|\ga_1|\,e^{-\be|\ga_1|}{1\over|\ga_1|}\sum_{\ga_2\not\sim\ga_1}|\ga_2|
    e^{-\be|\ga_2|}\dots {1\over|\ga_{k-1}|} \sum_{\ga_k\not\sim\ga_{k-1}}
    e^{-\be|\ga_{k}|}\nonumber\\
  &\le& e^{-t}\,\sum_{i=0}^\infty {t^i \over i!} \sum_{k\ge i}\alpha_0\,
  \alpha^{k-1}\nonumber\\
  &=& {\alpha_0\over \alpha(1-\alpha)}\, e^{-(1-\alpha)t}.  \label{802}
\end{eqnarray}

{\bf 2.} In Fern\'andez, Ferrari and Garcia (1998) we offered an alternative
proof of part (iii), based on a the computation of the exponential moment of
the total population of a subcritical single-type branching process, which
dominates the space width. However this proof works in a smaller range of
$\beta$.

\section{Proof of Theorem \ref{18}}\label{SS6}
 
The following theorem shows that the condition $\alpha<\infty$ implies
the hypothesis of Theorem \ref{20}.  This proves Theorem \ref{18}.

\begin{thm} \label{hat} If $\alpha <\infty$, then for all $x$ and positive $t$
  the set $A^{x,0}[-t,0]$ has a finite number of cylinders with probability
  one.
\end{thm}

\proof Let $C$ be a cylinder with basis $\ga$ and birth-time $0$. Recall the
definition of $U(C')$ just before display \reff{25} and define 
\begin{equation}
  \label{25b}
  \widetilde\psi^\ga_t(\th)= \Bigl|\Bigl\{C'\in \bar\B^C:
  \basis(C')=\th,-t< \birth(U(C'))\Bigr\}\Bigr|.
\end{equation}
The process $\widetilde\psi^\ga_t$ signals all contours born in $[0,t]$ in
the process $\psi^\ga_t$. Notice that for $x\in \ga$,
\begin{equation}
  \label{238}
 \Bigl|\A^{x,0}[-t,0]\Bigr|\;\le\;
\Bigl|\B^{x,0}[-t,0]\Bigr|\;=\;
\Bigl|\Bigl\{C\in \B^{x,0}: \birth(C)\in
  [-t,0]\Bigr\}\Bigr|\;\le\;
\sum_\th\widetilde\psi^\ga_t(\th)
\end{equation}
We prove that this is finite with probability one by showing it has a finite
mean. Indeed, reasoning as in the previous section,
\begin{equation}
  \label{25c}
  \label{25d}
\E\left(\sum_\th\widetilde\psi^\ga_t(\th)\right) \;=\;
\sum_\th \bigl(e^{\mm t}\bigr)(\ga,\th) \;\le\; |\ga|\,
  e^{t\alpha} \;<\; \infty
\end{equation}
if $\alpha<\infty$. \square

\section{Proof of Theorem \ref{18a}} \label{SS7}

We prove that the hypothesis of Theorem \ref{18a} imply those of Theorem
\ref{225}. The different parts of Theorem \ref{18a} follow by combining
Theorem \ref{225} and the space-width and time-length estimations of Theorem
\ref{thm:hall}. 

\paragraph{Existence and uniqueness}
In Theorem \ref{thm:hall} (i) the condition $\alpha(\beta)<1$ was
shown to imply lack of backward oriented percolation.  This, plus part
(1) of Theorem \ref{225}, proves (1) of Theorem \ref{18a}.

\paragraph{Exponential time convergence}
Inequality \reff{334} follows
from \reff{336}, \reff{46a} and \reff{335} by choosing $b=(2-\alpha)^{-1}$.
\ \square

\paragraph{Exponential space convergence}
In view of \reff{81a} it suffices to bound 
\begin{equation}
  \label{620}
  \P\Bigl( \sw(\A(\supp(f)))\, \ge\, d(\supp(f),\Lambda^c) \Bigr)\le
  \sum_{x\in\supp(f)}\P\Bigl( \sw(\A^{x,0})\, \ge\, d(\{x\},\Lambda^c)\Bigr). 
\end{equation}
By \reff{46c}, this is bounded by
\begin{equation}
  \label{a620}
   \sum_{x\in\supp(f)}
   \frac{\alpha_0(\widetilde\beta)}{1-\alpha(\widetilde\beta)}\, 
    e^{-(\beta-\widetilde\beta)\,d(\{x\},\Lambda^c)} .\ \square
\end{equation}

\paragraph{Exponential mixing}
We shall use part (4) of Theorem \ref{225}.  We first show, in the
next lemma, that $\alpha<\infty$ implies the existence of an $h$ such
that there is non-oriented percolation in the interval $(0,h)$.

\begin{lemma}
  \label{500}
  For all $h>0$ such that
  \begin{equation}
   \alpha\, h < 1,
  \end{equation}
  the probability that there is no (non-oriented) percolation in $(0,h)$ is
  one.
\end{lemma}

\proof Analogously to Section \ref{mbp}, we dominate the construction of the
set ancestors of the non oriented percolation process by a multitype branching
process.  In this branching process, the mean number of ancestors $\theta$ of
a contour $\gamma$ is
\begin{equation}
  \label{501}
 \overline\mm(\ga,\th) = h\,\one\{\th\not\sim\ga\}\,  e^{-\be|\th|}
\end{equation}
As in Lemma \ref{16k}, this branching process is subcritical if
$\overline\alpha = \alpha\,h<1$. \square
\bigskip

We prove now \reff{101}. From \reff{505} and \reff{506}:
\begin{equation}
  \label{63}
  \mu_\Lambda (fg) - \mu_\Lambda f\, \mu_\Lambda g
=\E(f(\eta_0)g(\eta_0))- \E f(\eta_0)\,\E g(\eta_0).
\end{equation}

By \reff{73b} it is enough to bound
\begin{equation}
  \label{507}
  \P\Bigl(\A(\supp(f)) \, \not\sim\, \widehat\A(\supp(g))\Bigr)
\end{equation}
where $\widehat\A(\supp(g))$ has the same distribution as $\A(\supp(g))$ but
is independent of $\A(\supp(f))$. This is bounded by
\begin{equation}
  \label{530}
  \sum_{\scriptstyle x\in \supp(f),\atop\scriptstyle y\in \supp(g)}
  \P\Bigl(\A^{x,0} \, \not\sim\, \widehat\A^{y,0}\Bigr) \le \sum_{\scriptstyle
    x\in \supp(f),\atop\scriptstyle y\in \supp(g)}
  \P\Bigl(\sw(\A^{x,0})+\sw(\widehat\A^{y,0})\ge |x-y|\Bigr)
\end{equation}
Using the following inequality which is valid for independent random variables
$S_1$ and $S_2$,
\begin{equation}
  \label{80}
  \P(S_1+S_2 \ge \ell) \le \sum_{j=1}^\ell \P(S_1\ge j)\,\P(S_2\ge \ell-j)
\end{equation}
and the exponential decay of \reff{46c} we get the decay stated in
\reff{101}. \square

\paragraph{Central Limit Theorem}
We apply the central limit theorem for stationary mixing random fields
proven by Bolthausen (1982). Let $X_x = \tau_x f$. Let ${\cal A}_\Lambda$ be the
sigma algebra generated by $\{X_x:x\in \Lambda\}$. Define
\begin{eqnarray}
  \label{100}
  \alpha_{k,\ell}(n) &=& \sup\Bigl\{|\P(A_1\cap A_2) -\P(A_1)\P(A_2)|: A_1 \in
  {\cal A}_{\Lambda_1}, A_2 \in
  {\cal A}_{\Lambda_2}, \nonumber\\
  &&\qquad\qquad\qquad\qquad|\Lambda_1|\le k\,,|\Lambda_2|\le \ell\,,\,
  d(\Lambda_1,\Lambda_2) \ge n\Bigr\}.
\end{eqnarray}

The simplified version of Bolthausen theorem stated in Remark 1 pag 1049 of
his paper says that if there exists a $\delta>0$ such that
$||X_x||_{2+\delta}<\infty$ and
\begin{equation}
  \label{93}
  \sum_{n=1}^\infty n^{d-1} (\alpha_{2,\infty}(n))^{\delta/(2+\delta)} <\infty
\end{equation}
then $D<\infty$ and \reff{90} holds. Hence, it suffices to show that
$\alpha_{2,\infty}(n)$ decays exponentially fast with $n$. We can write
\begin{equation}
  \label{94}
  \alpha_{2,\infty}(n) = \sup_{a, g_1, g_2} \,\left|\mu (g_1 g_2 ) - \mu g_1 \,  
    \mu g_2\right|
\end{equation}
where the supremum is taken over the set of $a\in\Z^d$, $g_1$ in the set of
indicator functions with support on $\supp(f) \cup \tau_a\supp(f)$ and $g_2$ in
the set of indicator functions with support in
\begin{equation}
  \label{111}
  \bigcup\,\Bigl\{\tau_y\supp(f): y \in\Z^d \hbox{ and } |y-x|\ge
  n\,\forall x\in \supp (f)\Bigr\}.
\end{equation}
By \reff{101}
\begin{eqnarray}
  \label{120}
  \alpha_{2,\infty}(n) &\le& 2\,(M_2)^2\,\sum_{\scriptstyle x\in
    \supp(g_1),\atop\scriptstyle y\in  \supp(g_2)}
  \,|x-y| e^{-M_3 |x-y|}\nonumber\\
&\le& 4\,(M_2)^2\,|\supp(f)|\,\sum_{|y|\ge n-2|\supp(f)|} e^{-M_3 |y|}
\end{eqnarray}
because $|\supp (g_1)|\le 2\,|\supp(f)|$ and $\|g_1\|_\infty=\|g_2\|_\infty=1$.
Hence, $\alpha_{2,\infty}(n)$ decreases exponentially fast with $n$.\ \square

\def\x{\underline x}

\section{Proof of Theorem \ref{191}. Poisson approximation}\label{S62}

We define first a common probability space where all processes
$\{\eta^\be:\be*<\be\le \infty\}$ can simultaneously be constructed.
For each $\ga\in\widetilde\G= \cup_j\widetilde\G_j$ let $\NN_\ga$ be a
marked Poisson process on $\R^{d+2}$ of rate $1$. The event points of
this process are denoted $(u,t,r,s)$, where $u\in \R^d$, $t,r\in \R$
and $s\in \R^+$. The coordinate $t$ is interpreted as time while the
coordinate $r$ is later used to tune the rate of the projected process
$(u,t)$. The coordinate $s$ ---the \emph{mark}--- is an exponential
random variable with mean one independent of everything (used later to
determine the lifetime of the corresponding point/cylinder).  Denote
by $\P$ the product measure generated by
$(\NN_\ga:\ga\in\widetilde\G)$, and by $\E$ the corresponding
expectation. We identify the random counting measure $\NN_\ga$ with
the corresponding discrete random subset of $R^{d+1}\times\R^+$.

Fix a contour length $j$ and an inverse temperature $\beta$.  By counting only
those points in $\NN_{\th}$ whose $r$ coordinate is in $[0,e^{-\be(|\th|-j)}]$
we generate the $(d+1)$-dimensional marked process
\begin{equation}
  \label{904}
  \NN_{\th,\be} \;=\; \hbox { marked Poisson process of rate }
  e^{-\be(|\th|-j)}. 
\end{equation}
The life of each point $(u,t,s)\in\NN_{\th,\be}\times\R^+$ is 
the interval $[t,t+s]$.

Define a family of marked point processes indexed by $\beta$ and 
$\tau_x\th$, $x\in\Z^d$, $\th\in\widetilde G$ for Borel sets $I\subset
\R$ by
\begin{equation}
  \label{900}
  N_{\tau_x\th,\be}(I)\; =\; \NN_{\th,\be} (O(x e^{-\be j/d},e^{-\be
    j/d}/2)\times I )
\end{equation}
where for $y=(y_1,\dots,y_d)\in \Z^d$, $O$ is the $d$-dimensional ``rectangle''
\begin{equation}
  \label{901}
  O(y,\rho) = [y_1-\rho, y_1+\rho]\times\dots\times[y_d-\rho, y_d+\rho]
\end{equation}
Since the volume of $O(x e^{-\be j/d},e^{-\be j/d}/2)$ is $e^{-\be j}$, and
$\NN_{\th,\be}$ has rate $e^{-\be(|\th|-j)}$, the resulting process
$N_{\tau_x\th,\be}(I)$ is a one dimensional marked Poisson process of rate
$e^{-\be|\th|}$. The marks are the independent exponentially distributed random
variables of mean one, inherited from $\NN_{\th,\be}$. The point of this
construction is that all these Poisson processes are constructed
simultaneously as a function of the original $d+2$-dimensional Poisson processes.
  
Now we use the processes $N_{\tau_x\th,\be}$ to perform the graphical
construction of Section \ref{4.1}.  We call $\C^\be$ be the family of
cylinders so obtained. Let $\xi^\be_t$ be the free network of Section
\ref{free}.  and $\eta^\be_t$ the loss networks of \reff{201}. As in
\reff{eq:muga} and \reff{505} these processes have invariant distributions
$\mu^0_\beta$ and $\mu_\beta$ respectively.

Let $V$ be a $d$-dimensional rectangle as in the statement of
the theorem. Let
\begin{equation}
  \label{630}
  M_{\ga,\beta}^0(V) = \sum_{x\in V\cdot e^{\be |\ga|/d}}
  \xi^\beta_0(\tau_x\ga) = \sum_{x\in V\cdot e^{\be
  |\ga|/d}}\,\sum_{C\in\C^\be} \one\{\basis(C)=\tau_x\ga, \life(C)\ni 0\} 
\end{equation}
as in \reff{201}. The super-label zero in the left hand side indicates we are
dealing with the free process $\xi^\beta_t$, while the sublabel zero in the
right hand side indicates time zero. The family
$(M^0_{\ga,\beta}(V):\ga\in\widetilde G_j)$ consists of $|\widetilde G_j|$
independent Poisson random variables with mean:
\begin{equation}
  \label{631}
  \E M^0_{\ga,\beta}(V) =  |V\cdot e^{\be |\ga|/d}|\,e^{-\be|\ga|}.
\end{equation}

By \reff{210} $\eta^\be_0(\ga)$ constructed with the cylinders in $\C^\be$ is
$\mu_\be$ distributed. Thus we can use $\eta^\be_0(\ga)$ in the definition
\reff{920} of $M_{\ga,\beta}$. By \reff{a623}
\begin{equation}
  \label{923}
\eta^\beta_0(\ga)\le\xi^\beta_0(\ga).
\end{equation}
 Hence, $ M_{\ga,\beta}(V) \le
M^0_{\ga,\beta}(V)$.

The joint construction also implies
\begin{eqnarray}
  \label{123}
  \P(M^0_{\ga,\beta}(V) - M_{\ga,\beta}(V) \ge 1) 
  &\le& \sum_{x\in V\cdot e^{\be |\ga|/d}}
  \P\Bigl(\xi^\beta_0(\tau_x\ga)-\eta^\beta_0(\tau_x\ga)\ge 1\Bigr) \nonumber\\
  &\le& \sum_{x\in V\cdot e^{\be |\ga|/d}}
  \biggl[\P\Bigl(\xi^\be_0(\tau_x\ga)\ge 1,\eta^\be_0(\tau_x\ga)=0
  \Bigr)\nonumber\\
  &&\qquad\qquad{} + \P\Bigl(\xi^\be_0(\tau_x\ga)\ge 2,\eta^\be_0(\tau_x\ga)=1
  \Bigr)\biggl].
\end{eqnarray}
From the construction, for any $\th \in \G$, 
\begin{equation}
  \label{634}
  \Bigl\{\xi^\be_0(\th)\ge 1,\eta^\be_0(\th)=0\Bigr\} \;\subset \;
  \Bigl\{\C^\beta: \C^\beta\ni C \hbox{ with } \basis(C)=\th\, \life(C)\ni 0,
  \A_1^C \ne \emptyset\Bigr\}.
\end{equation}
The probability of this last event is bounded by
\begin{equation}
  \label{123a}
  \P(\xi^\be_0(\th)\ge 1) \;\P\Bigl(\sum_{\th'} \b_1^{\th}(\th')\ge
  1\Bigr).
\end{equation}
Since $\P(\xi^\be_0(\th)\ge 1) = 1- \exp(e^{-\be |\th|})\leq e^{-\be |\th|}$,
the rhs of \reff{123a} is bounded above by
\begin{equation}
  \label{124}
  e^{-\be |\th|} \sum_{\th':\th'\not\sim\th} e^{-\be|\th'|}\;\le\;
  e^{-\be|\th|}\,|\th|\,\alpha(\beta).
\end{equation}
On the other hand, 
\begin{equation}
  \label{910}
  \P\Bigl(\xi^\be_0(\th)\ge 2,\eta^\be_0(\th)=1 \Bigr)\; \le
  \;\P\Bigl(\xi^\be_0(\th)\ge 2\Bigr)\; \le \; {1\over 2}\, e^{-2\be|\th|}
\end{equation}
From \reff{123}--\reff{910} we get
\begin{equation}
  \label{911}
 \P(M^0_{\ga,\beta}(V) - M_{\ga,\beta}(V) \ge 1)\;\le\; \Bigl|V\cdot e^{\be
    |\ga|/d}\Bigr|\,e^{-\be|\ga|}\,\Bigl(|\ga|\,\alpha(\beta)\,\,+\, {1\over
    2}\, e^{-\be|\ga|}\Bigr)\;\sim\;e^{-2d\be}.
\end{equation}

To finish the proof of \reff{125} we must show that $M^0_{\ga,\beta}$
is close to a Poisson process. For $|\ga|=j$, let $M_{\ga,\infty}$ count
those points of the Poisson process $\NN_{\ga,\infty}$ whose life contains the
origin. The process $M_{\ga,\infty}$ is a Poisson process in $\R^d$ of rate
one.  This is because the lifetimes are independent exponentials of mean one
and for $|\ga|=j$, $\NN_{\ga,\infty}$ is a Poisson process of rate one. The
family $\{M_{\ga,\infty}:\ga\in \widetilde G_j\}$ inherits independence from
$\{\NN_{\ga}:\ga\in\widetilde\G\}$.
For $J\subset \Z^d$ let
\begin{equation}
  \label{926}
  J:a\; = \;\{r\in \R^{d}: ra \in J+[1/2\,,1/2]^d\}\; \subset\; \R^d.
\end{equation}
By definition, 
\begin{equation}
  \label{758}
  M^0_{\ga,\beta}(V) \;=\; M_{\ga,\infty}\Bigl((V\cdot e^{\be j/d}):e^{\be
    j/d}\Bigr).
\end{equation}
Then, as $(V\cdot e^{\be j/d}):e^{\be j/d} \subset V$,
\begin{eqnarray}
  \label{925}
  \P\Bigl(M^0_{\ga,\beta}(V)- M_{\ga,\infty}(V) \neq 0\Bigr) &\le& \P\biggl(
  M_{\ga,\infty}\Bigl(V \setminus \bigl[(V\cdot e^{\be j/d}):e^{\be 
    j/d}\bigr]\Bigr) >
  0\biggr)\nonumber\\
  &\le& \Bigl|V \setminus \bigl[(V\cdot e^{\be j/d}):e^{\be
    j/d}\bigr]\Bigr|\nonumber\\
  &\le& 2d\, |V|^{(d-1)/d}\, e^{-\be j/d}.
\end{eqnarray}
 
Inequality \reff{125} follows from \reff{911} and \reff{925}.

Proposition I.2 of Neveu (1977) and the comments below the statement of the
Proposition say that the distribution on finite unions of $d$- dimensional
finite-volume rectangles is enough to characterize a point process.  Since the
estimates \reff{125} can be easily extended to finite unions of rectangles,
the weak convergence follows. \ \square

\baselineskip 12pt

\parskip 1mm 
\parindent 1cm 
\section*{Acknowledgments} We thank Jean Bricmont, Joel L. Lebowitz, Enzo
Olivieri, Errico Presutti, Timo Sepp\"al\"ainen, Alan Sokal and
Bernard Ycart for instructive discussions and criticism.

This work was partially supported by FAPESP 95/0790-1 (Projeto Tem\'atico
``Fen\^omenos cr\'\i ticos e processos evolutivos e sistemas em equil\'\i
brio'') CNPq, FINEP (N\'ucleo de Excel\^encia ``Fen\^omenos cr\'\i ticos em
probabilidade e processos estoc\'asticos'' PRONEX-177/96).

\parskip 0pt

\vskip 6mm
\obeylines
Roberto Fern\'andez, \hfill Pablo A. Ferrari
IEA USP, \hfill IME USP, 
Av. Prof. Luciano Gualberto, \hfill Caixa Postal 66281, 
Travessa J, 374 T\'erreo \hfill 05389-970 - S\~{a}o Paulo,
05508-900 - S\~{a}o Paulo,\hfill BRAZIL
BRAZIL \hfill email: {\tt pablo@ime.usp.br}
email: {\tt rf@ime.usp.br} \hfill http://www.ime.usp.br/\~{}pablo
\vskip 5mm

Nancy L. Garcia
IMECC, UNICAMP, Caixa Postal 6065, 
13081-970 - Campinas SP 
BRAZIL
email: {\tt nancy@ime.unicamp.br}
http://www.ime.unicamp.br/\~{}nancy

\end{document}